\documentclass[reqno]{amsart}
\usepackage{amssymb}
\usepackage{amsmath}
\usepackage{amsfonts}
\usepackage{graphicx}
\usepackage{subfig}
\usepackage{eurosym}
\usepackage{amssymb}
\usepackage{amsmath}
\usepackage{amsfonts}
\usepackage{graphicx}
\usepackage{subfig}

\newtheorem{theorem}{Theorem}
\theoremstyle{plain}

\newtheorem{definition}{Definition}

\newtheorem{lemma}{Lemma}

\newtheorem{remark}{Remark}

\numberwithin{equation}{section}

\begin{document}
\title[Batman University]{Smarandache curves and their properties on null
curves in lightlike cone space $%
\mathbb{Q}
_{2}^{3}$}
\author{Fatma ALMAZ}
\address{department of mathematics, faculty of arts and sciences, batman
university, batman/ t\"{u}rk\.{ı}ye orcid: 0000-0002-1060-7813}
\email{fatma.almaz@batman.edu.tr}
\author{B\.{ı}lal TEKYOLDA\c{S}}
\address{department of mathematics, faculty of arts and sciences, batman
university, batman/ t\"{u}rk\.{ı}ye, orcid: 0009-0003-2630-4846}
\email{bilaltekyoldas63@gmail.com}
\subjclass{53B30, 53C50, 14H45}
\keywords{Lightlike cone $Q_{3}^{2}\subset E_{2}^{4}$, natural Frenet frame,
smarandache curves.}
\thanks{This paper is in final form and no version of it will be submitted
for publication elsewhere.}

\begin{abstract}

This study addresses the differential geometric representation and intrinsic kinematics of Smarandache curves parametrized along lightlike trajectories within the specific degenerate hypersurface geometry of the lightlike cone space $\mathbb{Q}_{2}^{3} \subset E_{2}^{4}$. Driven by mathematical physics motivations—specifically string world-sheet dynamics and massless particle trajectories in general relativity the indefinite metric tensor of signature $(2,2)$ is leveraged to construct a more complex and rich geometric framework than those found in standard Euclidean or low-dimensional Minkowski spaces. Utilizing the unique invariant structure of the natural null Frenet frame $\{x, \xi, N, W\}$, we systematically formulate the bending, torsion, and higher-order structural variations for seven distinct configurations of Smarandache curves, generated via advanced linear and transcendental combinations of the tangent, transversal normal, and binormal vector fields. Beyond the straightforward derivation of new curvature functions, this paper analytically characterizes the global geometric transitions of these curves. We explicitly establish the critical degeneration criteria under which a Smarandache curve collapses into a singular state or preserves lightlike behavior, and we define the exact constraints necessary for these trajectories to form generalized helices. These structural insights bridge the gap between constrained pseudo-Riemannian submanifolds and the evolving geometric modelling of theoretical particle kinematics.

\end{abstract}

\maketitle

\section{Introduction}

Differential geometry stands as a cornerstone of modern mathematical research, providing powerful intrinsic and extrinsic toolsets for analyzing the structural topographies of curves and surfaces. Beyond its role in pure mathematical frameworks, this field serves as the geometric language for prominent physical disciplines, including general relativity, topological fluid dynamics, and string theory. While classical Riemannian geometry standardizes our intuitions through a positive-definite metric tensor, the transition to semi-Riemannian manifolds offers a much larger and radically complex geometric landscape by relaxing this requirement \cite{14}. On these pseudo-Riemannian manifolds, the inner product allows vectors to be spacelike, timelike, or null depending on their zero or non-zero causal norms. This causal classification reveals rich, degenerate, and directional geometric phenomena that deviate entirely from standard Euclidean intuitions. 

The present study focuses specifically on a four-dimensional, two-index semi-Riemannian space $E_{2}^{4}$ and its embedded degenerate hypersurface, namely the lightlike cone $\mathbb{Q}_{2}^{3}$. Metric spaces with such high-index null vector structures possess distinct invariant characters that separate them fundamentally from the classical Minkowski space-time $E_{1}^{3}$ or $E_{1}^{4}$. Consequently, these spaces have emerged as highly active domains of investigation within both mathematical physics (e.g., super gravity configurations and localized target spaces in string theory) and modern submanifold theory \cite{9, 10}. Null curves within this framework, whose tangent vectors maintain a vanishing causal norm while remaining distinct from the zero vector, act as the exact mathematical representations for the space-time world-lines of massless particles or photons travelling at the absolute speed of light \cite{12, 13}.

In recent years, the geometric generation of new trajectories from a parent curve via the linear combination of its movement frame vectors commonly referred to as Smarandache curves has attracted significant attention. This architectural concept allows researchers to analyze the high-order structural and acceleration invariants of an underlying curve from completely novel perspectives. To map these behaviors, pioneering researchers have described diverse frameworks for special Smarandache curves across different settings. For instance, foundational definitions and theorems were established within classical Euclidean spaces by Ali \cite{3}, and subsequently extended to higher-dimensional spaces and alternative moving frames by Akan and Yayl\i\ \cite{2}, as well as \"{O}nder and Bilgi\c{c} \cite{16}. In the context of indefinite metric tensors, structural formulations were advanced within Minkowski space-times by Turgut and Y{\i}lmaz \cite{19}, \c{C}etin et al. \cite{8}, \"{O}nder and \"{O}zcan \cite{17}, and Y{\i}lmaz \cite{20}. Furthermore, Abdel-Aziz and Saad examined these curves under non-flat restrictions within Galilean space geometries \cite{1}, while Ali established explicit models for timelike Smarandache configurations derived from spacelike helices \cite{4}.

The structural behavior of curves restricted to degenerate hypersurfaces requires specialized mathematical treatment due to the vanishing nature of the induced metric tensor. To address this, structural characterizations of Smarandache curves within null cone spaces of varying dimensions have been investigated in the literature. Notably, Almaz and K\"{u}lahc{\i} systematically formulated special Smarandache frameworks directly inside the null cone \cite{5}, explored their approximation boundaries \cite{6}, and advanced the geometric assessment criteria of their invariants within the specific stage $Q^2$ \cite{11}. Similarly, \"{O}zt\"{u}rk et al. investigated the exact behavior of Smarandache trajectories lying on the lightlike cone within Minkowski 3-space \cite{15}. This line of geometric inquiry builds upon the fundamental axiom combinations introduced in the structural study of Smarandache Geometries by Ashbacher \cite{7}, which highlights how combining non-Euclidean axiom sets gives rise to non-degenerate structural behaviors.

Despite these collective insights, a significant gap remains in the literature regarding the global kinematics and structural stability of Smarandache curves derived from null curves inside the higher-dimensional lightlike cone space $\mathbb{Q}_{2}^{3} \subset E_{2}^{4}$. The indefinite metric structure of signature $(2,2)$ produces a unique geometric setting where the cross-interaction between tangent vectors and transversal null normals behaves differently than in low-index Lorentz-Minkowski spaces \cite{18}. 

To bridge this gap, the main objective of this paper is to conduct a comprehensive kinematical analysis of null curves embedded in $\mathbb{Q}_{2}^{3}$ by employing the natural null Frenet frame $\{x, \xi, N, W\}$. Building upon this foundational framework, we systematically derive the bending, torsion, and high-order geometric invariants of seven distinct, advanced Smarandache curve configurations (e.g., combinations involving tangent, transversal normal, and binormal vectors). Moving beyond routine algebraic formulations, we provide a complete analysis of their structural transitions by explicitly evaluating the critical parameters under which these newly generated trajectories maintain lightlike properties or encounter singular degenerations. This comprehensive review aims to expand the body of knowledge in modern curve theory within indefinite metric manifolds and to establish a unified geometric foundation for theoretical physics applications.

\section{Preliminaries}

In this section, we establish the rigorous geometric conventions, indefinite metric structures, and the operational moving frame equations necessary for characterizing the intrinsic kinematics of null trajectories within degenerate hypersurfaces.

The $4-$dimensional flat pseudo-Euclidean space $E_{2}^{4}$ is defined by equipping the real vector space $\mathbb{R}^4$ with a standard non-degenerate, indefinite metric tensor of signature $(2,2)$. For any two vectors $x=(x_{1},x_{2},x_{3},x_{4})$ and $y=(y_{1},y_{2},y_{3},y_{4})$ in $E_{2}^{4}$, the symmetric bilinear form is explicitly given by
\begin{equation}
\left\langle x,y\right\rangle =\sum_{i=1}^{2} x_{i}y_{i}-\sum_{j=3}^{4} x_{j}y_{j}. \label{eq:metric}
\end{equation}

Let $M$ be a smooth submanifold embedded in $E_{2}^{4}$. The causal character of $M$ is governed by the restriction of the ambient metric form. Specifically, $M$ is classified as a timelike, spacelike, or degenerate (lightlike) submanifold if the induced metric tensor is pseudo-Riemannian, Riemannian, or a degenerate quadratic form on the tangent space $T_pM$ at each point $p \in M$, respectively \cite{9}.

Let $c \in E_{2}^{4}$ denote a fixed point acting as the origin, and let $r>0$ be a real constant. The pseudo-Riemannian lightlike cone hypersurface with the center point $c$ is defined as the set
\begin{equation}
\mathbb{Q}_{2}^{3}(c,r)=\{ \beta \in E_{2}^{4} : \left\langle \beta - c, \beta - c \right\rangle = 0 \}, \label{eq:cone}
\end{equation}
which constitutes a canonical degenerate hypersurface in $E_{2}^{4}$. The tangent space $T_p\mathbb{Q}_{2}^{3}$ at any point contains a radical subspace $\text{Rad}(T_p\mathbb{Q}_{2}^{3}) = T_p\mathbb{Q}_{2}^{3} \cap (T_p\mathbb{Q}_{2}^{3})^\perp \neq \{0\}$, identifying it as a structurally degenerate pseudo-Riemannian submanifold.

A smooth parametrized curve $\gamma : I \longrightarrow \mathbb{Q}_{2}^{3} \subset E_{2}^{4}$, defined on an open interval $I \subset \mathbb{R}$, is said to be a null curve if its velocity vector field $\gamma^{\prime}(t) = \xi(t) = \frac{d\gamma(t)}{dt}$ satisfies the null condition at each point, meaning $\left\langle \xi(t), \xi(t) \right\rangle = 0$ for all $t \in I$ while $\xi(t) \neq 0$. 

To rule out non-regular structural collapse along the trajectory, we assume that the null curve $\gamma$ is parametrized by its natural pseudo-arclength parameter $s$, such that the acceleration vector field satisfies the normalizing structural constraint \cite{10}:
\begin{equation}
\left\langle \gamma^{\prime\prime}(s), \gamma^{\prime\prime}(s) \right\rangle = 1. \label{eq:pseudo_arclength}
\end{equation}

For any null curve $\gamma$ restricted to the lightlike cone $\mathbb{Q}_{2}^{3}$, there exists a unique, non-degenerate moving frame field $\{\gamma(s), \xi(s), N(s), W(s)\}$ along $\gamma(s)$, known as the natural null Frenet frame field \cite{10}. The constitutive structural equations governed by the indefinite inner product tensor require these frame vectors to satisfy the following rigid pseudo-orthogonality relations:
\begin{equation}
\begin{aligned}
\left\langle \gamma, \gamma \right\rangle &= 0, \quad \left\langle \gamma, \xi \right\rangle = 0, \quad \left\langle \gamma, N \right\rangle = 0, \quad \left\langle \gamma, W \right\rangle = 1, \\
\left\langle \xi, \xi \right\rangle &= 0, \quad \left\langle \xi, N \right\rangle = 1, \quad \left\langle \xi, W \right\rangle = 0, \\
\left\langle N, N \right\rangle &= 0, \quad \left\langle N, W \right\rangle = 0, \quad \left\langle W, W \right\rangle = 0,
\end{aligned} \label{eq:frenet_orthogonal}
\end{equation}
where $N(s)$ represents the unique lightlike transversal vector field associated with the tangent direction $\xi(s)$, and $W(s)$ is the unique lightlike transversal vector field orthogonal to the position vector $\gamma(s)$. 

Consequently, the differential geometry and intrinsic kinematics of the null curve $\gamma$ on the degenerate hypersurface $\mathbb{Q}_{2}^{3}$ are completely determined by the following system of natural null Frenet equations:
\begin{equation}
\begin{aligned}
\gamma^{\prime}(s) &= \xi(s), \\
\xi^{\prime}(s) &= h(s)\xi(s) + \kappa_{1}(s)\gamma(s), \\
N^{\prime}(s) &= -h(s)N(s) + \kappa_{2}(s)\gamma(s) - W(s), \\
W^{\prime}(s) &= -\kappa_{2}(s)\xi(s) - \kappa_{1}(s)N(s),
\end{aligned} \label{eq:frenet_eqs}
\end{equation}
where the invariant functional coefficients $h(s), \kappa_{1}(s),$ and $\kappa_{2}(s)$ constitute the fundamental curvature functions of the null curve $\gamma$ on $\mathbb{Q}_{2}^{3}$. These geometric invariants are explicitly determined via the following projections:
\begin{equation}
h(s) = \left\langle \xi^{\prime}(s), N(s) \right\rangle, \quad \kappa_{1}(s) = \left\langle \xi^{\prime}(s), W(s) \right\rangle, \quad \kappa_{2}(s) = \left\langle N^{\prime}(s), W(s) \right\rangle. \label{eq:curvatures}
\end{equation}

\begin{definition} \cite{9, 18}
Let $E_{2}^{4}$ be the flat pseudo-Euclidean space endowed with the symmetric bilinear form of signature $(2,2)$ as defined in \eqref{eq:metric}. For any lightlike vector field $V=(v_{1},v_{2},v_{3},v_{4}) \in E_{2}^{4}$, there exists a canonical isotropic transformation induced by the algebraic product of the split-quaternion structure of the ambient space. The natural orthogonal null companion vector fields to $V$, denoted by $V_{13}^{\bot}$ and $V_{14}^{\bot}$ respectively, are characterized as the isotropic dual configurations determined via the explicit component permutations:
\begin{equation}
V_{13}^{\bot}=(v_{2},-v_{1},v_{4},-v_{3}) \quad \text{and} \quad V_{14}^{\bot}=(v_{2},-v_{1},-v_{4},v_{3}). \label{eq:companion_permutation}
\end{equation}

Geometrically, these fields represent the fundamental isotropic reflections that map the localized lightlike direction $V$ into its complementary totally degenerate planes, satisfying the strict boundary validation condition $\left\langle V, V_{1i}^{\bot} \right\rangle = 0$ for $i=3,4$. This algebraic operator guarantees the preservation of the structural index pairing across the split-signature geometry of the manifold.
\end{definition}

\begin{lemma} \cite{10}
Let $U$ and $V$ be two arbitrary null vector fields defined on the split-signature pseudo-Euclidean manifold $E_{2}^{4}$. The structural isotropic isometry induced by the natural companion operators $V_{1i}^{\bot}$ satisfies the following bilinear inner product invariants under the fixed permutation index $i \in \{3,4\}$:
\begin{itemize}
    \item[\rm i)] $\left\langle U, V \right\rangle = \left\langle U_{1i}^{\bot}, V_{1i}^{\bot} \right\rangle$,
    \item[\rm ii)] $\left\langle U, V_{1i}^{\bot} \right\rangle = -\left\langle U_{1i}^{\bot}, V \right\rangle$.
\end{itemize}

This isometric invariance guarantees that the companion projection maps smoothly across the totally degenerate linear subspaces of the space.
\end{lemma}

By synthesizing these structural properties, we establish a rigorous mathematical foundation for the underlying null trajectories embedded within the lightlike cone space $\mathbb{Q}_{2}^{3}$, facilitating the explicit formulation and stability analysis of the Smarandache curve classes.

\section{Representation of Smarandache Curves in the Lightlike Cone $\mathbb{Q}_{2}^{3}$}

In this section, null curves defined according to the natural Frenet frame within the lightlike cone $\mathbb{Q}_{2}^{3}$ in space, inspired by Smarandache curves, are examined. Furthermore, the curvature functions $h, \kappa _{1}$, and $\kappa _{2}$ are given for different Smarandache curve configurations.

Let $\gamma =(\gamma _{1},\gamma _{2},\gamma _{3},\gamma_{4}):I\longrightarrow \mathbb{Q}_{2}^{3}\subset E_{2}^{4}$ be a null parent trajectory embedded within the lightlike cone. By leveraging the isotropic definition of the cone geometry, the position vector satisfies $\left\langle \gamma ,\gamma \right\rangle =0$, which explicitly expands to the quadratic constraint $-\gamma _{1}^{2}-\gamma _{2}^{2}+\gamma _{3}^{2}+\gamma _{4}^{2} = 0$. Rearranging this equation into the symmetric ratio form $\gamma _{1}^{2}-\gamma_{3}^{2}=\gamma _{4}^{2}-\gamma _{2}^{2}$, we introduce the operational functional parametrizations $m$ and $n$ via the following rigid relational system:
\begin{equation}
\frac{\gamma _{1}+\gamma _{3}}{\gamma _{4}+\gamma _{2}}=-\frac{\gamma_{2}-\gamma _{4}}{\gamma _{1}-\gamma _{3}}=m, \tag{3.1}
\end{equation}
\begin{equation}
\frac{\gamma _{1}-\gamma _{3}}{\gamma _{2}-\gamma _{4}}=-\frac{\gamma_{2}+\gamma _{4}}{\gamma _{1}+\gamma _{3}}=n \tag{3.2}
\end{equation}
and the spatial scaling parameter functions correspond to:
\begin{equation}
\gamma _{2}+\gamma _{4}=h=2g. \tag{3.3}
\end{equation}

Hence, by performing back-substitution via the simultaneous systems governed by (3.1), (3.2), and (3.3), one can write the following explicit coordinate coupling equations:
\begin{equation}
\gamma _{1}+\gamma _{3}=-\frac{h}{n}=-\frac{2g}{n}, \tag{3.4a}
\end{equation}
\begin{equation}
\gamma _{1}-\gamma _{3}=hm=2gm; \quad \gamma _{2}-\gamma _{4}=-hm^{2}=-2gm^{2}. \tag{3.4b}
\end{equation}

Extracting individual localized coordinate expressions directly from (3.4a) and (3.4b), one obtains the intermediate parametrizations:
\begin{equation*}
\gamma _{1}=\frac{g}{n}(mn-1); \quad \gamma _{2}=gm(1-m); \quad \gamma _{3}=-\frac{g}{n}(mn+1); \quad \gamma _{4}=gm(1+m).
\end{equation*}

Hence, by setting the continuous function $n=-\frac{g}{f}$, the entire position vector field of the null parent curve $\gamma $ can be written elegantly in terms of the component profiles $f$ and $g$ as follows:
\begin{equation}
\gamma =(f+mg, \, g-mf, \, f-mg, \, g+mf), \tag{3.5}
\label{eq:gamma_profile_final_3}
\end{equation}
which intrinsically guarantees the balance relation:
\begin{equation}
\gamma _{1}^{2}+\gamma _{2}^{2}=\gamma _{3}^{2}+\gamma _{4}^{2}. \tag{3.6}
\end{equation}

Since $\gamma : I \rightarrow \mathbb{Q}_{2}^{3}$ is a regular lightlike trajectory, its velocity vector field satisfies the null validation condition $\left\langle \gamma ^{\prime },\gamma ^{\prime }\right\rangle =0$, which directly implies the derivative norm equality:
\begin{equation}
\gamma _{1}^{\prime 2}+\gamma _{2}^{\prime 2}=\gamma _{3}^{\prime 2}+\gamma _{4}^{\prime 2}. \tag{3.7}
\end{equation}

Then, taking the algebraic product of the structural conditions presented in (3.6) and (3.7) yields:
\begin{equation*}
\left( \gamma _{1}^{2}+\gamma _{2}^{2}\right) \left( \gamma _{1}^{\prime 2}+\gamma _{2}^{\prime 2}\right) =\left( \gamma _{3}^{2}+\gamma_{4}^{2}\right) \left( \gamma _{3}^{\prime 2}+\gamma _{4}^{\prime 2}\right),
\end{equation*}
which expands through direct polynomial multiplication to the following form:
\begin{equation}
\gamma _{1}^{2}\gamma _{1}^{\prime 2}+\gamma _{1}^{2}\gamma _{2}^{\prime 2}+\gamma _{2}^{2}\gamma _{1}^{\prime 2}+\gamma _{2}^{2}\gamma _{2}^{\prime 2}=\gamma _{3}^{2}\gamma _{3}^{\prime 2}+\gamma _{3}^{2}\gamma _{4}^{\prime 2}+\gamma _{4}^{2}\gamma _{3}^{\prime 2}+\gamma _{4}^{2}\gamma _{4}^{\prime 2}. \tag{3.8}
\end{equation}

Furthermore, the rigid pseudo-orthogonality requirement of the manifold requires $\left\langle \gamma ,\gamma ^{\prime }\right\rangle =0$, generating the invariant relation $\left( \gamma _{1}\gamma _{1}^{\prime }+\gamma _{2}\gamma _{2}^{\prime}\right) ^{2}=\left( \gamma _{3}\gamma _{3}^{\prime }+\gamma _{4}\gamma_{4}^{\prime }\right) ^{2}$. Expanding this square structures the system as:
\begin{equation}
\left( \gamma _{1}\gamma _{1}^{\prime }\right) ^{2}+\left( \gamma _{2}\gamma _{2}^{\prime }\right) ^{2}+2\gamma _{1}\gamma _{1}^{\prime }\gamma _{2}\gamma _{2}^{\prime }=\left( \gamma _{3}\gamma _{3}^{\prime }\right) ^{2}+\left( \gamma _{4}\gamma _{4}^{\prime }\right) ^{2}+2\gamma _{3}\gamma _{3}^{\prime }\gamma _{4}\gamma _{4}^{\prime }. \tag{3.9}
\end{equation}

From (3.8) and (3.9), one has
\begin{equation}
\gamma _{1}^{2}\gamma _{2}^{\prime 2}+\gamma _{2}^{2}\gamma _{1}^{\prime 2}-2\gamma _{1}\gamma _{1}^{\prime }\gamma _{2}\gamma _{2}^{\prime }=\gamma_{4}^{2}\gamma _{3}^{\prime 2}+\gamma _{3}^{2}\gamma _{4}^{\prime 2}-2\gamma_{3}\gamma _{3}^{\prime }\gamma _{4}\gamma _{4}^{\prime }. \tag{3.10}
\end{equation}

On the other hand, substituting the final parametric coordinates from (3.5) directly into the differential constraint (3.10) collapses the entire high-dimensional trajectory system into the following compact scalar form:
\begin{equation*}
(fg^{\prime }-f^{\prime }g)(1+m^{2})(f^{2}+g^{2})m^{\prime }=0.
\end{equation*}

Hence, assuming non-trivial functional components such that $f,g\neq 0,$ one can write the localized parameter constraints as:
\begin{equation}
fg^{\prime }-f^{\prime }g=\Omega \neq 0, \quad m^{\prime }=0, \tag{3.11}
\end{equation}
which explicitly proves that $m$ is a strict real constant along the entire null parent trajectory.

To analyse specific geometric configurations without loss of generality, we examine explicit profile constraints under parametrized coordinate mappings in the subsequent sub proofs.

\begin{definition}

Let $\gamma :I\rightarrow \mathbb{Q}_{2}^{3}\subset E_{2}^{4}$ be a null curve with non-zero curvatures $h^{1},\kappa _{1}^{1}$ and $\kappa _{2}^{1}$ and $\{\gamma (t),\xi (t),N(t),W(t)\}$ be moving frame on it. Then, $\gamma W-$Smarandache curve is defined with
\begin{equation}
\gamma _{\gamma W}=\sinh \psi _{1}\overrightarrow{\gamma }+\cosh \psi _{1}\overrightarrow{W}.  \tag{3.12}
\label{eq:gammaW_def_clean}
\end{equation}
\end{definition}

\begin{theorem}
Let $\gamma :I\rightarrow \mathbb{Q}_{2}^{3}\subset E_{2}^{4}$ is a null curve in $\mathbb{Q}_{2}^{3}$. Then, $h^{(1)},\kappa _{1}^{(1)}$ and $\kappa _{2}^{(1)}$ curvatures of the $\gamma W-$Smarandache curve are as follows, respectively
\begin{equation*}
\begin{split}
h^{(1)}&=\frac{-1}{(1+m^{2})\psi _{1}^{\prime }} \Big[ -b_{1}^{1}\left( \cosh \psi _{1}-m\sinh \psi _{1}\right) +b_{2}^{1}\left( \sinh \psi _{1}+m\cosh \psi _{1}\right) \\ 
&\quad -b_{3}^{1}\left( \cosh \psi _{1}+m\sinh \psi _{1}\right) +b_{4}^{1}\left( \sinh \psi _{1}-m\cosh \psi _{1}\right) \Big],
\end{split}
\end{equation*}
\begin{equation*}
\begin{split}
\kappa _{1}^{(1)}&=\frac{1}{(1+m^{2})} \Big[ -b_{1}^{1}\left( \sinh \psi _{1}-m\cosh \psi _{1}\right) +b_{2}^{1}\left( \cosh \psi _{1}+m\sinh \psi _{1}\right) \\ 
&\quad -b_{3}^{1}\left( \sinh \psi _{1}+m\cosh \psi _{1}\right) +b_{4}^{1}\left( \cosh \psi _{1}-m\sinh \psi _{1}\right) \Big],
\end{split}
\end{equation*}
\begin{equation*}
\begin{split}
\kappa _{2}^{(1)}&=\frac{-1}{(1+m^{2})} \Big[ -n_{1}^{1}\left( \sinh \psi _{1}-m\cosh \psi _{1}\right) +n_{2}^{1}\left( \cosh \psi _{1}+m\sinh \psi _{1}\right) \\ 
&\quad -n_{3}^{1}\left( \sinh \psi _{1}+m\cosh \psi _{1}\right) +n_{4}^{1}\left( \cosh \psi _{1}-m\sinh \psi _{1}\right) \Big],
\end{split}
\end{equation*}
where $M_{1}=\sqrt{-\psi _{1}^{\prime 2}+(\kappa _{1}^{2}-\kappa _{2}^{2})\cosh ^{2}\psi _{1}-\sinh ^{2}\psi _{1}+\kappa _{2}\sinh(2\psi _{1})}$,
\begin{equation*}
a_{1}^{1} =\frac{\psi _{1}^{\prime }\cosh \psi _{1}}{M_{1}}, \quad a_{2}^{1}=\frac{\sinh \psi _{1}-\kappa _{2}\cosh \psi _{1}}{M_{1}}, \quad a_{3}^{1} =\frac{-\kappa _{1}\cosh \psi _{1}}{M_{1}}, \quad a_{4}^{1}=\frac{\psi _{1}^{\prime }\sinh \psi _{1}}{M_{1}},
\end{equation*}
\begin{equation*}
\begin{aligned}
b_{1}^{1} &=\frac{a_{1}^{1\prime }+a_{2}^{1}\kappa _{1}+\kappa _{2}a_{3}^{1}}{M_{1}}, & b_{2}^{1}&=\frac{a_{2}^{1\prime }+a_{1}^{1}+a_{2}^{1}h-\kappa _{2}a_{4}^{1}}{M_{1}}, \\
b_{3}^{1} &=\frac{a_{3}^{1\prime }-a_{3}^{1}h-\kappa _{1}a_{4}^{1}}{M_{1}}, & b_{4}^{1}&=\frac{-a_{3}^{1}+a_{4}^{1\prime }}{M_{1}},
\end{aligned}
\end{equation*}
\begin{equation*}
\footnotesize
\begin{aligned}
n_{1}^{1} &=\frac{-1}{1+m^{2}}\left( \frac{\psi _{1}^{\prime \prime }}{\psi _{1}^{\prime 2}}\left( \cosh \psi _{1}-m\sinh \psi _{1}\right) +\sinh \psi _{1}-m\cosh \psi _{1}\right), \\
n_{2}^{1} &=\frac{1}{1+m^{2}}\left( \frac{\psi _{1}^{\prime \prime }}{\psi _{1}^{\prime 2}}\left( \sinh \psi _{1}+m\cosh \psi _{1}\right) +\cosh \psi _{1}+m\sinh \psi _{1}\right), \\
n_{3}^{1} &=\frac{1}{1+m^{2}}\left( \frac{\psi _{1}^{\prime \prime }}{\psi _{1}^{\prime 2}}\left( \cosh \psi _{1}+m\sinh \psi _{1}\right) +\sinh \psi _{1}+m\cosh \psi _{1}\right), \\
n_{4}^{1} &=\frac{-1}{1+m^{2}}\left( \frac{\psi _{1}^{\prime \prime }}{\psi _{1}^{\prime 2}}\left( \sinh \psi _{1}-m\cosh \psi _{1}\right) +\cosh \psi _{1}-m\sinh \psi _{1}\right).
\end{aligned}
\end{equation*}
\end{theorem}

\begin{proof}

Let $\gamma $ be a unit speed regular $\gamma W-$Smarandache curve as in \eqref{eq:gammaW_def_clean}. If one takes the derivative of the Smarandache curve according to arc length parameter and by using (2.4), one has
\begin{equation}
\begin{split}
\frac{d\gamma _{\gamma W}}{ds_{\gamma }}\frac{ds_{\gamma }}{ds}&=\xi _{\gamma _{\gamma W}}\frac{ds_{\gamma }}{ds}=\psi _{1}^{\prime }\cosh \psi _{1}\overrightarrow{\gamma } +(\sinh \psi _{1}-\kappa _{2}\cosh \psi _{1})\overrightarrow{\xi } \\
&\quad -\kappa _{1}\cosh \psi _{1}\overrightarrow{N}+\psi _{1}^{\prime }\sinh \psi _{1}\overrightarrow{W},
\end{split} \tag{3.13}
\label{eq:proof_deriv1_fixed}
\end{equation}
by taking the norm of (3.13), one gets
\begin{equation}
\frac{ds_{\gamma }}{ds}=\sqrt{\vert-\psi _{1}^{\prime 2}+(\kappa _{1}^{2}-\kappa _{2}^{2})\cosh ^{2}\psi _{1}-\sinh ^{2}\psi _{1}+\kappa _{2}\sinh(2\psi _{1})\vert}=M_{1}, \tag{3.14}
\label{eq:proof_M1_fixed}
\end{equation}
and the tangent vector of $\gamma _{\gamma W}$ is given as
\begin{equation}
\xi _{\gamma _{\gamma W}}=a_{1}^{1}\overrightarrow{\gamma }+a_{2}^{1}\overrightarrow{\xi }+a_{3}^{1}\overrightarrow{N}+a_{4}^{1}\overrightarrow{W}. \tag{3.15}
\label{eq:proof_tangent_fixed}
\end{equation}

By taking derivative (3.15), one can see that
\begin{equation}
\gamma _{\gamma W}^{\prime \prime }=\xi _{\gamma _{\gamma W}}^{\prime }=(b_{1}^{1}\overrightarrow{\gamma }+b_{2}^{1}\overrightarrow{\xi }+b_{3}^{1}\overrightarrow{N}+b_{4}^{1}\overrightarrow{W}), \tag{3.16}
\label{eq:proof_second_deriv_fixed}
\end{equation}

Considering the conditions $\gamma_1+\gamma_3\neq 0,$ $\gamma_2+\gamma_4\neq 0,$ $\gamma_2-\gamma_4\neq 0$ for the null curve and 
considering the special structural case where the coordinate functions are modeled as $\gamma_1=\sinh \psi_1; \gamma_4=\sinh \psi_1; f = \frac{\sinh \psi_1}{2}; g = \frac{\cosh \psi_1}{2}$, one gets $fg^{\prime}-f^{\prime}g = \Omega_1 = \frac{-\psi_1^{\prime}}{4}$. Also, from (3.11) $m = \text{constant}$ and by using Definition 1 and Lemma 1, one has
\begin{equation}
\gamma _{13}^{\bot }=(g-mf,-f-mg,g+mf,-f+mg), \tag{3.17}
\end{equation}
\begin{equation}
\gamma _{14}^{\bot }=(g-mf,-f-mg,-g-mf,f-mg) \tag{3.18}
\end{equation}
\begin{equation}
\gamma ^{\prime }=(f^{\prime }+mg^{\prime },g^{\prime }-mf^{\prime },f^{\prime }-mg^{\prime },g^{\prime }+mf^{\prime }). \tag{3.19}
\end{equation}

Now, let's construct the vector fields $N_{\gamma W}$ and $W_{\gamma W}$ \ 
\begin{equation}
N_{\gamma W}=\frac{1}{\left\langle \gamma _{14}^{\bot },\gamma ^{\prime }\right\rangle }\gamma _{14}^{\bot };W_{\gamma W}=\frac{1}{\left\langle \gamma _{14}^{\bot },\gamma ^{\prime }\right\rangle }\gamma _{14}^{\prime \bot }, \tag{3.20}
\end{equation}
also from equations (3.18) and (3.19) one has
\begin{equation}
\left\langle \gamma _{14}^{\bot },\gamma ^{\prime }\right\rangle =2(1+m^{2})\Omega _{1}, \tag{3.21}
\end{equation}
where $\Omega _{1}=\frac{-\psi _{1}^{\prime }}{4}$. Hence, from (3.20), (3.21) and $f=\frac{\sinh \psi _{1}}{2},g=\frac{\cosh \psi _{1}}{2}$, we get
\begin{equation}
N_{\gamma W}=\frac{-1}{(1+m^{2})\psi _{1}^{\prime }}\left( \begin{aligned} &\cosh \psi _{1}-m\sinh \psi _{1}, -\sinh \psi _{1}-m\cosh \psi _{1}, \\ &-\cosh \psi _{1}-m\sinh \psi _{1}, \sinh \psi _{1}-m\cosh \psi _{1} \end{aligned} \right), \tag{3.22}
\label{eq:N_explicit_fixed}
\end{equation}
\begin{equation}
W_{\gamma W}=\frac{1}{(1+m^{2})}\left( \begin{aligned} &\sinh \psi _{1}-m\cosh \psi _{1}, -\cosh \psi _{1}-m\sinh \psi _{1}, \\ &-\sinh \psi _{1}-m\cosh \psi _{1}, \cosh \psi _{1}-m\sinh \psi _{1} \end{aligned} \right). \tag{3.23}
\label{eq:W_explicit_fixed}
\end{equation}

By using (2.6) and direct calculations, one can obtain the curvature functions.
For example, the curvature functions for a null curve $ \gamma $. Hence, from(3.16)and (3.22),one gets the curvatures.
\end{proof}

\begin{figure*}[!htbp]
\centering
\includegraphics[width=0.95\textwidth]{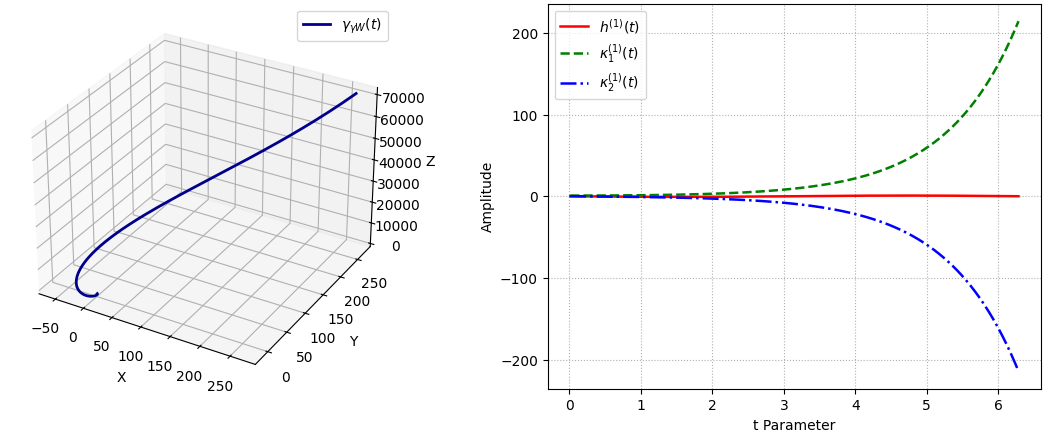}
\caption{The 3D spatial trajectory of the $\gamma W$-Smarandache curve and its  curvatures ($h^{(1)}(t)$, $\kappa_1^{(1)}(t)$, $\kappa_2^{(1)}(t)$)}
\label{fig:gammaW_trajectory_curvatures}
\end{figure*}

\begin{definition}

Let $\gamma :I\rightarrow \mathbb{Q}_{2}^{3}\subset E_{2}^{4}$ is a null curve with non-zero curvatures $h^{2},\kappa _{1}^{2}$ and $\kappa _{2}^{2}$ and $\{\gamma (t),\xi (t),N(t),W(t)\}$ be moving frame on it. Then, $\xi N-$Smarandache curve is defined with
\begin{equation}
\gamma _{\xi N}=\sinh \psi _{2}\overrightarrow{\xi }+\cosh \psi _{2}\overrightarrow{N}.  \tag{3.24}
\label{eq:xiN_def_clean}
\end{equation}
\end{definition}

\begin{theorem}

Let $\gamma :I\rightarrow \mathbb{Q}_{2}^{3}\subset E_{2}^{4}$ is a null curve in $\mathbb{Q}_{2}^{3}$. Then, $h^{(2)},\kappa _{1}^{(2)}$ and $\kappa _{2}^{(2)}$ curvatures of the $\xi N-$Smarandache curve are as follows, respectively
\begin{equation}
\begin{split}
h^{(2)}&=\frac{1}{(1+m^{2})\psi_2^{\prime}} \Big[ -n_{1}^{2}\left( \sinh \psi_{2}-m\cosh \psi_{2}\right) +n_{2}^{2}\left( \cosh \psi_{2}+m\sinh \psi_{2}\right) \\ 
&\quad -n_{3}^{2}\left( \cosh \psi_{2}+m\sinh \psi_{2}\right) +n_{4}^{2}\left( \sinh \psi_{2}-m\cosh \psi_{2}\right) \Big],
\end{split} \tag{3.24a}
\end{equation}
\begin{equation}
\begin{split}
\kappa_1^{(2)}&=\frac{1}{(1+m^{2})} \Big[ n_{1}^{2}\left( \cosh \psi_{2}-m\sinh \psi_{2}\right) +n_{2}^{2}\left( -\sinh \psi_{2}-m\cosh \psi_{2}\right) \\ 
&\quad -n_{3}^{2}\left( -\cosh \psi_{2}-m\cosh \psi_{2}\right) -n_{4}^{2}\left( \sinh \psi_{2}-m\cosh \psi_{2}\right) \Big],
\end{split} \tag{3.24b}
\end{equation}
\begin{equation}
\begin{split}
\kappa_2^{(2)}&=\frac{-1}{(1+m^{2})} \Big[ -\Upsilon_1^{2}\left( \cosh \psi_{2}-m\sinh \psi_{2}\right) +\Upsilon_2^{2}\left( \sinh \psi_{2}+m\cosh \psi_{2}\right) \\ 
&\quad -\Upsilon_3^{2}\left( \cosh \psi_{2}+m\sinh \psi_{2}\right) +\Upsilon_4^{2}\left( \sinh \psi_{2}-m\cosh \psi_{2}\right) \Big],
\end{split} \tag{3.24c}
\end{equation}
where the underlying spatial parametrization metric factor satisfies the split relation:
\begin{equation}
\begin{split}
M_{2}^{2} &= -\kappa_{1}^{2}\sinh^{2}\psi_{2} + (1-\kappa_{2}^{2})\cosh^{2}\psi_{2} - \psi_{2}^{\prime 2} + h^2 \\
&\quad - \sinh(2\psi_{2})(\kappa_{1}\kappa_{2} + 2\psi_{2}^{\prime}h),
\end{split} \tag{3.24d}
\end{equation}
such that the scaling invariant factor resolves precisely to $M_{2}=\sqrt{M_{2}^{2}}$, and the structural vector contraction arrays $a_i^2$ expand through:
\begin{equation}
\begin{aligned}
a_{1}^{2} &=\kappa_{1}\sinh \psi_{2}+\kappa_{2}\cosh \psi_{2}, \\
a_{2}^{2} &=\psi_{2}^{\prime}\cosh \psi_{2}+h\sinh \psi_{2}, \\
a_{3}^{2} &=\psi_{2}^{\prime}\sinh \psi_{2}-h\cosh \psi_{2},
\end{aligned} \tag{3.24e}
\end{equation}
while the operational normal frame components $n_i^2$ yield the simultaneous equations:
\begin{equation}
\begin{aligned}
n_{1}^{2} &=\frac{1}{M_{2}^{2}}\left( -\frac{M_{2}^{\prime }}{M_{2}}a_{1}^{2}+a_{1}^{2\prime }+a_{2}^{2}\kappa_{1}+a_{3}^{2}\kappa_{2}\right), \\
n_{2}^{2} &=-\frac{M_{2}^{\prime }}{M_{2}^{3}}a_{2}^{2}+\frac{1}{M_{2}^{2}}\left( a_{1}^{2}+a_{2}^{2\prime }+ha_{2}^{2}+\kappa_{2}\cosh \psi_{2}\right), \\
n_{3}^{2} &=-\frac{M_{2}^{\prime }}{M_{2}^{3}}a_{3}^{2}+\frac{1}{M_{2}^{2}}\left( a_{3}^{2\prime }-ha_{3}^{2}+\kappa_{1}\cosh \psi _{2}\right), \\
n_{4}^{2} &=\frac{-1}{M_{2}^{2}}\left( \psi_{2}^{\prime}\sinh \psi_{2}+a_{3}^{2}\right),
\end{aligned} \tag{3.24f}
\end{equation}
and the derivative vector fields $\Upsilon_{i}^{2}$ tracking the higher-order variations of the lightlike frame satisfy the direct relations:
\begin{equation}
\begin{aligned}
\Upsilon_{1}^{2} &= \frac{d}{dt}\left[ \frac{\sinh \psi_{2} - m\cosh \psi_{2}}{(1+m^{2})\psi_{2}^{\prime}} \right], &
\Upsilon_{2}^{2} &= \frac{d}{dt}\left[ \frac{- \cosh \psi_{2} - m\sinh \psi_{2}}{(1+m^{2})\psi_{2}^{\prime}} \right], \\
\Upsilon_{3}^{2} &= \frac{d}{dt}\left[ \frac{-\sinh \psi_{2} - m\cosh \psi_{2}}{(1+m^{2})\psi_{2}^{\prime}} \right], &
\Upsilon_{4}^{2} &= \frac{d}{dt}\left[ \frac{\cosh \psi_{2} - m\sinh \psi_{2}}{(1+m^{2})\psi_{2}^{\prime}} \right].
\end{aligned} \tag{3.24g}
\end{equation}
\end{theorem}

\begin{proof}

Let $\gamma $ be a unit speed regular $\xi N-$Smarandache curve as in \eqref{eq:xiN_def_clean}. If we take the derivative of the Smarandache curve according to arc length parameter and from (2.4), one has
\begin{equation}
\begin{split}
\frac{d\gamma _{\xi N}}{ds_{\gamma }}\frac{ds_{\gamma }}{ds}&=\xi _{\gamma _{\xi N}}\frac{ds_{\gamma }}{ds}=(\kappa _{1}\sinh \psi _{2}+\kappa _{2}\cosh \psi _{2})\overrightarrow{\gamma } \\
&\quad +(\psi _{2}^{\prime }\cosh \psi _{2}+h\sinh \psi _{2})\overrightarrow{\xi } +(\psi _{2}^{\prime }\sinh \psi _{2}-h\cosh \psi _{2})\overrightarrow{N}-\cosh \psi _{2}\overrightarrow{W},
\end{split} \tag{3.25}
\label{eq:proof_xiN_deriv_fixed}
\end{equation}
by taking the norm of (3.25), one finds
\begin{equation}
\frac{ds_{\gamma }}{ds}=\sqrt{-\kappa _{1}^{2}\sinh ^{2}\psi _{2}+(1-\kappa _{2}^{2})\cosh ^{2}\psi _{2}-\psi _{2}^{\prime 2}+h^{2}-\sinh(2\psi _{2})(\kappa _{1}\kappa _{2}+2\psi _{2}^{\prime }h)}=M_{2}, \tag{3.26}
\label{eq:proof_xiN_M2_fixed}
\end{equation}
and the tangent vector of $\gamma _{\xi N}$ is
\begin{equation}
\xi _{\gamma _{\xi N}}=\frac{1}{M_{2}}\left( a_{1}^{2}\overrightarrow{\gamma }+a_{2}^{2}\overrightarrow{\xi }+a_{3}^{2}\overrightarrow{N}-\cosh \psi _{2}\overrightarrow{W}\right). \tag{3.27}
\label{eq:proof_xiN_tangent_fixed}
\end{equation}

By taking derivative (3.27), one can write that
\begin{equation}
\gamma _{\xi N}^{\prime \prime }=\xi _{\gamma _{\xi N}}^{\prime }=(n_{1}^{2}\overrightarrow{\gamma }+n_{2}^{2}\overrightarrow{\xi }+n_{3}^{2}\overrightarrow{N}+n_{4}^{2}\overrightarrow{W}), \tag{3.28}
\end{equation}

Considering the conditions for the null curve and from equations $f=\frac{\cosh \psi _{2}}{2},g=\frac{\sinh \psi _{2}}{2}$, one gets 
\begin{equation}
fg^{\prime }-f^{\prime }g=\Omega _{2}=\frac{\psi _{2}^{\prime }}{4}. \tag{3.29}
\end{equation} 

Also, from Definition 1 and Lemma 1 under $m = \text{const.}$, the transversal frame fields expand dynamically:
\begin{equation}
N_{\xi N}=\frac{1}{(1+m^{2})\psi _{2}^{\prime }}\left( \begin{aligned} &\sinh \psi _{2}-m\cosh \psi _{2},-\cosh \psi _{2}-m\sinh \psi _{2}, \\ &-\sinh \psi _{2}-m\cosh \psi _{2},\cosh \psi _{2}-m\sinh \psi _{2} \end{aligned} \right) \tag{3.30}
\label{eq:N_xiN_explicit_fixed}
\end{equation}
\begin{equation}
W_{\xi N}=\frac{-1}{(1+m^{2})}\left( \begin{aligned} &\cosh \psi _{2}-m\sinh \psi _{2},-\sinh \psi _{2}-m\cosh \psi _{2}, \\ &-\cosh \psi _{2}-m\sinh \psi _{2},\sinh \psi _{2}-m\cosh \psi _{2} \end{aligned} \right) \tag{3.31}
\end{equation}
and from (2.6) the curvature functions $ \kappa_{1}^{2} $ and $ h^{2} $ for a null curve $ \gamma $ are obtained and by taking derivative of (3.30), one gets \eqref{eq:N_xiN_explicit_fixed} isolates $N_{\xi N}^{\prime }=\left( \Upsilon _{1}^{2},\Upsilon _{2}^{2},\Upsilon _{3}^{2},\Upsilon _{4}^{2}\right)$. Then, the curvature $ \kappa_{2}^{2} $ is proven easily.
\end{proof}

\begin{figure*}[!htbp]
\centering
\includegraphics[width=0.95\textwidth]{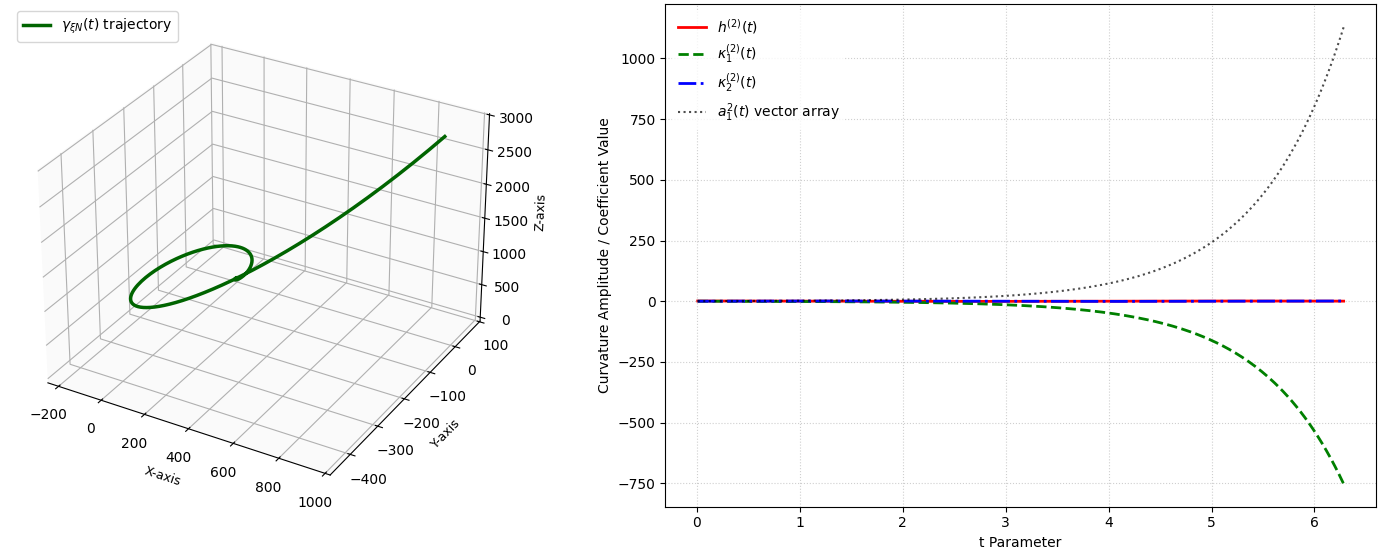}
\caption{The 3D trajectory and the curvatures of the $\xi N$-Smarandache curve.}
\label{fig:xiN_trajectory_curvatures}
\end{figure*}

\begin{definition}
Let $\gamma :I\rightarrow \mathbb{Q}_{2}^{3}\subset E_{2}^{4}$ is a null curve with non-zero curvatures $h^{3},\kappa _{1}^{3}$ and $\kappa _{2}^{3}$ and $\{\gamma (t),\xi (t),N(t),W(t)\}$ be moving frame on it. Then, $WN-$Smarandache curve is defined with
\begin{equation}
\gamma _{WN}=\sin \psi _{3}\overrightarrow{N}+\cos \psi _{3}\overrightarrow{W}.  \tag{3.32}
\label{eq:WN_def_clean}
\end{equation}
\end{definition}

\begin{theorem}
Let $\gamma :I\rightarrow \mathbb{Q}_{2}^{3}\subset E_{2}^{4}$ is a null curve in $\mathbb{Q}_{2}^{3}$. Then, $h^{(3)},\kappa _{1}^{(3)}$ and $\kappa _{2}^{(3)}$ curvatures of the $WN-$Smarandache curve are as follows, respectively
\begin{equation}
\begin{split}
h^{(3)}&=\frac{1}{(1+m^{2})\psi _{3}^{\prime }} \Big[ a_{1}^{3}\left( \cos \psi _{3}-m\sin \psi _{3}\right) -a_{2}^{3}\left( \sin \psi _{3}+m\cos \psi _{3}\right) \\ 
&\quad +a_{3}^{3}\left( \cos \psi _{3}+m\sin \psi _{3}\right) -a_{4}^{3}\left( \sin \psi _{3}-m\cos \psi _{3}\right) \Big],
\end{split} \tag{3.32a}
\end{equation}
\begin{equation}
\begin{split}
\kappa_1^{(3)}&=\frac{1}{(1+m^{2})} \Big[ a_{1}^{3}\left( \sin \psi _{3}+m\cos \psi _{3}\right) -a_{2}^{3}\left( -\cos \psi _{3}+m\sin \psi _{3}\right) \\ 
&\quad +a_{3}^{3}\left( \sin \psi _{3}-m\cos \psi _{3}\right) +a_{4}^{3}\left( \cos \psi _{3}+m\sin \psi _{3}\right) \Big],
\end{split} \tag{3.32b}
\end{equation}
\begin{equation}
\begin{split}
\kappa_2^{(3)}&=\frac{1}{(1+m^{2})} \Big[ n_{1}^{3}\left( \sin \psi _{3}+m\cos \psi _{3}\right) -n_{2}^{3}\left( -\cos \psi _{3}+m\sin \psi _{3}\right) \\ 
&\quad +n_{3}^{3}\left( \sin \psi _{3}-m\cos \psi _{3}\right) +n_{4}^{3}\left( \cos \psi _{3}+m\sin \psi _{3}\right) \Big],
\end{split} \tag{3.32c}
\end{equation}
where the underlying spatial parametrization metric factor satisfies the split relation:
\begin{equation}
\begin{split}
M_{3}&=\Big( -\kappa _{2}^{2}+\sin ^{2}\psi _{3}\left[h^{2}+(1+\psi _{3}^{\prime })^{2}\right] +\cos ^{2}\psi _{3}(\psi _{3}^{\prime }-\kappa _{1})^{2} \\
&\quad -h(\psi _{3}^{\prime }-\kappa _{1})\sin(2\psi _{3})\Big)^{\frac{1}{2}},
\end{split} \tag{3.32d}
\end{equation}
and the structural components $a_i^3$ representing the tangent projections adapt precisely to the following localized system:
\begin{equation}
\begin{aligned}
a_{1}^{3} &= \frac{-\kappa_{1}\sin \psi_{3} + (\psi_{3}^{\prime} - \kappa_{1})\cos \psi_{3}}{M_{3}}, \\
a_{2}^{3} &= \frac{h\sin \psi_{3} + (1 + \psi_{3}^{\prime})\cos \psi_{3}}{M_{3}}, \\
a_{3}^{3} &= \frac{\kappa_{2}\sin \psi_{3} - h\cos \psi_{3}}{M_{3}}, \\
a_{4}^{3} &= \frac{\sin \psi_{3}}{M_{3}},
\end{aligned} \tag{3.32e}
\end{equation}
while the operational secondary invariants $n_i^3$ yield the simultaneous differential constraints:
\begin{equation}
\begin{aligned}
n_{1}^{3} &= \frac{1}{M_{3}^{2}}\left( -\frac{M_{3}^{\prime }}{M_{3}}a_{1}^{3} + a_{1}^{3\prime } + a_{2}^{3}\kappa_{1} + a_{3}^{3}\kappa_{2} \right), \\
n_{2}^{3} &= -\frac{M_{3}^{\prime }}{M_{3}^{3}}a_{2}^{3} + \frac{1}{M_{3}^{2}}\left( a_{1}^{3} + a_{2}^{3\prime } + h a_{2}^{3} + \kappa_{2}\cos \psi_{3} \right), \\
n_{3}^{3} &= -\frac{M_{3}^{\prime }}{M_{3}^{3}}a_{3}^{3} + \frac{1}{M_{3}^{2}}\left( a_{3}^{3\prime } - h a_{3}^{3} + \kappa_{1}\cos \psi_{3} \right), \\
n_{4}^{3} &= \frac{-1}{M_{3}^{2}}\left( \psi_{3}^{\prime}\sin \psi_{3} + a_{3}^{3} \right).
\end{aligned} \tag{3.32f}
\end{equation}
\end{theorem}

\begin{proof}
Let $\gamma $ be a unit speed regular $WN-$Smarandache curve as in \eqref{eq:WN_def_clean}. If we take the derivative of the Smarandache curve according to arc length parameter and by using (2.4), one has
\begin{equation}
\frac{d\gamma _{WN}}{ds_{\gamma }}\frac{ds_{\gamma }}{ds} = \xi_{WN}\frac{ds_{\gamma }}{ds} = b_{1}^{3}\overrightarrow{\gamma }+b_{2}^{3}\overrightarrow{\xi }+b_{3}^{3}\overrightarrow{N}+b_{4}^{3}\overrightarrow{W}, \tag{3.33}
\label{eq:proof_WN_deriv_fixed}
\end{equation}
where 
\begin{eqnarray*}
b_{1}^{3} &=&\kappa _{2}\sin \psi _{3};b_{2}^{3}=-\kappa _{2}\cos \psi _{3};
\\
b_{3}^{3} &=&(\psi _{3}^{\prime }-\kappa _{1})\cos \psi _{3}-h\sin \psi
_{3};b_{4}^{3}=-(1+\psi _{3}^{\prime })\sin \psi _{3}.
\end{eqnarray*}
and by taking the norm of (3.33), one gets
\begin{equation}
\begin{split}
\frac{ds_{\gamma }}{ds} &= \Big(\vert -\kappa _{2}^{2}+\sin ^{2}\psi _{3}\left[h^{2}+(1+\psi _{3}^{\prime })^{2}\right] +\cos ^{2}\psi _{3}(\psi _{3}^{\prime }-\kappa _{1})^{2} \\
&\quad -h(\psi _{3}^{\prime }-\kappa _{1})\sin(2\psi _{3}) \Big\vert)^{\frac{1}{2}} = M_{3},
\end{split} \tag{3.34}
\label{eq:proof_WN_M3_fixed}
\end{equation}
and the tangent vector of $\gamma _{WN}$ is given as
\begin{equation}
\xi _{WN}=\frac{1}{M_{3}}\left( b_{1}^{3}\overrightarrow{\gamma }+b_{2}^{3}\overrightarrow{\xi }+b_{3}^{3}\overrightarrow{N}+b_{4}^{3}\overrightarrow{W}\right). \tag{3.35}
\label{eq:proof_WN_tangent_fixed}
\end{equation}

By taking derivative (3.35), one can obtain that
\begin{equation}
\gamma _{WN}^{\prime \prime }=\xi _{\gamma _{WN}}^{\prime }=a_{1}^{3}\overrightarrow{\gamma }+a_{2}^{3}\overrightarrow{\xi }+a_{3}^{3}\overrightarrow{N}+a_{4}^{3}\overrightarrow{W}, \tag{3.36}
\end{equation}
where $  a_{i}^{3}=\frac{1}{M_{3}} b_{i}^{3}$. Considering conditions for the null curve and choosing profiles $f=\frac{\sin \psi _{3}}{2},g=\frac{\cos \psi _{3}}{2}$, one gets
\begin{equation}
fg^{\prime }-f^{\prime }g=\Omega _{3}=\frac{-\psi _{3}^{\prime }}{4}. \tag{3.37}
\end{equation}

Under $m = \text{const.}$, from Definition 1 and by using (3.18), (3.19), (3.20), (3.37), one gets 
\begin{equation}
N_{WN}=\frac{-1}{(1+m^{2})\psi _{3}^{\prime }}\left( \begin{aligned} &\cos \psi _{3}-m\sin \psi _{3}, \, -\sin \psi _{3}-m\cos \psi _{3}, \\ &-\cos \psi _{3}-m\sin \psi _{3}, \, \sin \psi _{3}-m\cos \psi _{3} \end{aligned} \right) \tag{3.38}
\label{eq:N_WN_explicit_fixed}
\end{equation}
\begin{equation}
W_{WN}=\frac{1}{(1+m^{2})}\left( \begin{aligned} &-\sin \psi _{3}-m\cos \psi _{3}, \, -\cos \psi _{3}+m\sin \psi _{3}, \\ &\sin \psi _{3}-m\cos \psi _{3}, \, \cos \psi _{3}+m\sin \psi _{3} \end{aligned} \right). \tag{3.39}
\end{equation}

Executing the projections from \eqref{eq:curvatures} completes the verification for $h^{(3)}, \kappa_{1}^{(3)},$ and $\kappa_{2}^{(3)}$.
\end{proof}
\begin{figure*}[!htbp]
\centering
\includegraphics[width=0.95\textwidth]{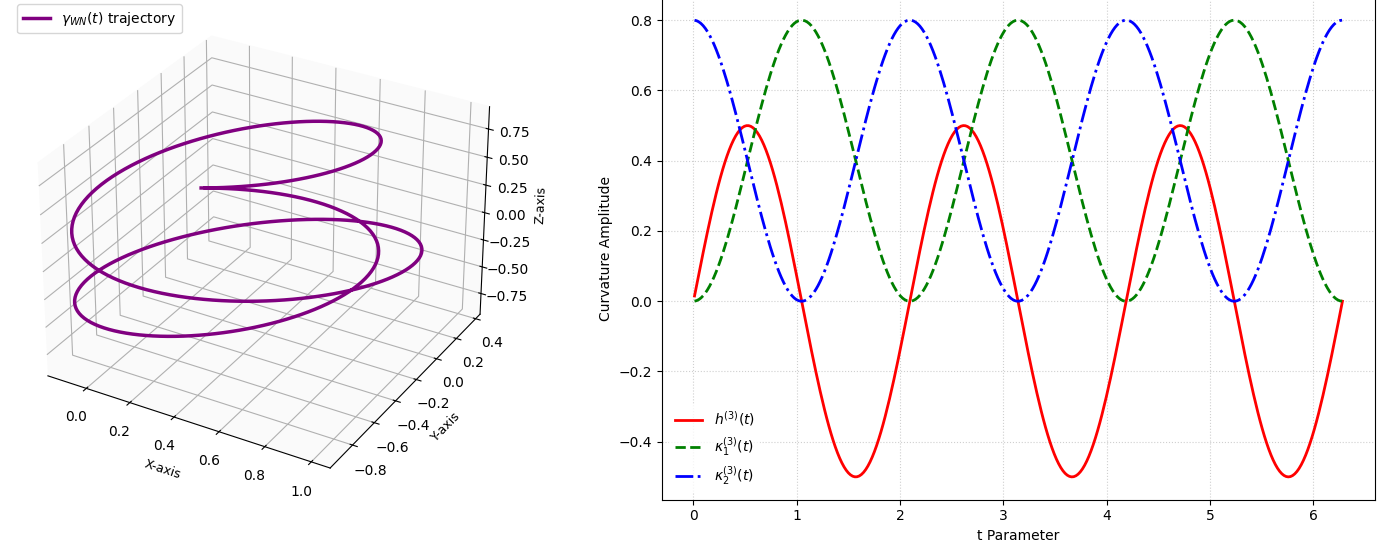}
\caption{The 3D trajectory and the curvatures of the $WN$-Smarandache curve.}
\label{fig:WN_trajectory_curvatures}
\end{figure*}

\begin{definition}
Let $\gamma :I\rightarrow \mathbb{Q}_{2}^{3}\subset E_{2}^{4}$ is a null curve with non-zero curvatures $h^{4},\kappa _{1}^{4}$ and $\kappa _{2}^{4}$ and $\{\gamma (t),\xi (t),N(t),W(t)\}$ be moving frame on it. Then, $\gamma \xi N-$Smarandache curve is defined with
\begin{equation}
\gamma _{\gamma \xi N}=\frac{1}{\sinh \Phi _{2}}\left(\sin \Phi _{1}\overrightarrow{\gamma }+\cos \Phi _{1}\overrightarrow{\xi }+\cosh \Phi _{2}\overrightarrow{N}\right).  \tag{3.40}
\label{eq:gamma_xiN_def_clean}
\end{equation}
\end{definition}

\begin{theorem}
Let $\gamma :I\rightarrow \mathbb{Q}_{2}^{3}\subset E_{2}^{4}$ is a null curve in $\mathbb{Q}_{2}^{3}$. Then, $h^{(4)},\kappa _{1}^{(4)}$ and $\kappa _{2}^{(4)}$ curvatures of the $\gamma \xi N-$Smarandache curve are as follows, respectively
\begin{equation}
\begin{split}
h^{(4)} &= \frac{1}{4(1+m^{2})\Omega _{4}\sinh \Phi _{2}} \Big[ -b_{1}^{4}\left( \cos \Phi _{1}-m\left( \sin \Phi _{1}+\cosh \Phi_{2}\right) \right) \\ 
&\quad +b_{2}^{4}\left( \sin \Phi _{1}+\cosh \Phi _{2}+m\cos \Phi _{1}\right) -b_{3}^{4}\left( \cos \Phi _{1}+m\left( \sin \Phi _{1}+\cosh \Phi_{2}\right) \right) \\ 
&\quad +b_{4}^{4}\left( \sin \Phi _{1}+\cosh \Phi _{2}-m\cos \Phi _{1}\right) \Big],
\end{split} \tag{3.40a}
\end{equation}
\begin{equation}
\begin{split}
\kappa_{1}^{(4)} &= \frac{-1}{2(1+m^{2})\Omega _{4}} \Big[ -b_{1}^{4}\left( h_{1}^{4\prime }+h_{2}^{4}\kappa_{1}+h_{3}^{4}\kappa_{2}\right) -b_{2}^{4}\left( h_{2}^{4\prime }+h_{1}^{4}+h_{2}^{4}h-h_{4}^{4}\kappa_{2}\right) \\ 
&\quad +b_{3}^{4}\left( h_{3}^{4\prime }-h_{3}^{4}h-h_{4}^{4}\kappa_{1}\right) +b_{4}^{4}\left( h_{4}^{4\prime }-h_{3}^{4}\right) \Big],
\end{split} \tag{3.40b}
\end{equation}
\begin{equation}
\begin{split}
\kappa_{2}^{(4)} &= \frac{1}{2(1+m^{2})\Omega _{4} } \Big[ -\left( n_{1}^{4\prime }+n_{2}^{4}\kappa _{1}+n_{3}^{4}\kappa _{2}\right) \left( h_{1}^{4\prime }+h_{2}^{4}\kappa _{1}+h_{3}^{4}\kappa _{2}\right) \\ 
&\quad -\left( n_{2}^{4\prime }+n_{1}^{4}+n_{2}^{4}h-n_{4}^{4}\kappa _{2}\right) \left( h_{2}^{4\prime }+h_{1}^{4}+h_{2}^{4}h-h_{4}^{4}\kappa _{2}\right) \\ 
&\quad +\left( n_{3}^{4\prime }-n_{3}^{4}h-n_{4}^{4}\kappa _{1}\right) \left( h_{3}^{4\prime }-h_{3}^{4}h-h_{4}^{4}\kappa _{1}\right) + \left( n_{4}^{4\prime }-n_{3}^{4}\right) \left( h_{4}^{4\prime }-h_{3}^{4}\right) \Big],
\end{split} \tag{3.40c}
\end{equation}
where the underlying parametrization metric factor satisfies the split relation:
\begin{equation}
\begin{split}
M_4 &= \Big( -(\gamma _{1}^{\prime }+\kappa _{1}\gamma _{2})^{2}-(\gamma_{1}+\gamma _{2}^{\prime }+h\gamma _{2}+\kappa _{2}\gamma _{3})^{2} \\
&\quad +(\gamma_{3}^{\prime }-h\gamma _{3})^{2}+\gamma _{3}^{2} \Big)^{\frac{1}{2}},
\end{split} \tag{3.40d}
\end{equation}
and the orthogonal companion normal components $h_i^4$ expand through the symmetric arrangement:
\begin{equation}
\begin{aligned}
h_{1}^{4} &= \frac{\cos \Phi_{1} - m(\sin \Phi_{1} + \cosh \Phi_{2})}{2\sinh \Phi_{2}},h_{2}^{4}= \frac{-(\sin \Phi_{1} + \cosh \Phi_{2}) - m\cos \Phi_{1}}{2\sinh \Phi_{2}}, \\
h_{3}^{4}&= \frac{-\cos \Phi_{1} - m(\sin \Phi_{1} + \cosh \Phi_{2})}{2\sinh \Phi_{2}},h_{4}^{4}= \frac{(\sin \Phi_{1} + \cosh \Phi_{2}) - m\cos \Phi_{1}}{2\sinh \Phi_{2}},
\end{aligned} \tag{3.40e}
\end{equation}
while the operational movement framework quantities $n_i^4$ resolve into the final differential matrix blocks:
\begin{equation}
\begin{aligned}
n_{1}^{4} &= \frac{h_{1}^{4}}{2(1+m^2)\Omega_4}, &
n_{2}^{4} &= \frac{h_{2}^{4}}{2(1+m^2)\Omega_4}, \\
n_{3}^{4} &= \frac{h_{3}^{4}}{2(1+m^2)\Omega_4}, &
n_{4}^{4} &= \frac{h_{4}^{4}}{2(1+m^2)\Omega_4}.
\end{aligned} \tag{3.40f}
\end{equation}
\end{theorem}

\begin{proof}
Let $\gamma_{\gamma \xi N}$ be defined as in (3.40). Differentiating with respect to the parameter $s$, the velocity vector field satisfies
\begin{equation}
\begin{split}
\frac{d\gamma _{\gamma \xi N}}{ds_{\gamma }}\frac{ds_{\gamma }}{ds} &=\xi_{\gamma \xi N}\frac{ds_{\gamma }}{ds}=(\gamma _{1}^{\prime }+\kappa_{1}\gamma _{2})\overrightarrow{\gamma } \\
&\quad +(\gamma _{1}+\gamma _{2}^{\prime }+h\gamma _{2}+\kappa _{2}\gamma _{3})\overrightarrow{\xi } \\
&\quad +(\gamma _{3}^{\prime }-h\gamma _{3})\overrightarrow{N}-\gamma _{3}\overrightarrow{W}.
\end{split} \tag{3.41}
\end{equation}

Evaluating the norm yields the parametrization factor:
\begin{equation}
\begin{split}
\frac{ds_{\gamma }}{ds} &= \Big( -(\gamma _{1}^{\prime }+\kappa _{1}\gamma _{2})^{2}-(\gamma_{1}+\gamma _{2}^{\prime }+h\gamma _{2}+\kappa _{2}\gamma _{3})^{2} \\
&\quad +(\gamma_{3}^{\prime }-h\gamma _{3})^{2}+\gamma _{3}^{2} \Big)^{\frac{1}{2}} = M_{4},
\end{split} \tag{3.42}
\end{equation}
and the localized unit tangent vector field maps into the explicit form:
\begin{equation}
\xi _{\gamma \xi N}=\frac{1}{M_{4}}\left[ (\gamma _{1}^{\prime }+\kappa _{1}\gamma _{2})\overrightarrow{\gamma }+(\gamma _{1}+\gamma _{2}^{\prime }+h\gamma _{2}+\kappa _{2}\gamma _{3})\overrightarrow{\xi } +(\gamma _{3}^{\prime }-h\gamma _{3})\overrightarrow{N}-\gamma _{3}\overrightarrow{W} \right]. \tag{3.43}
\end{equation}

Differentiating (3.43) yields the variation vector field:
\begin{equation}
\gamma _{\gamma \xi N}^{\prime \prime }=\xi _{\gamma _{\gamma \xi N}}^{\prime }=b_{1}^{4}\overrightarrow{\gamma }+b_{2}^{4}\overrightarrow{\xi }+b_{3}^{4}\overrightarrow{N}+b_{4}^{4}\overrightarrow{W}, \tag{3.44}
\end{equation}

Setting the coordinate forms $f=\frac{\sin \Phi _{1}+\cosh \Phi _{2}}{2\sinh \Phi _{2}}$ and $g=\frac{\cos \Phi _{1}}{2\sinh \Phi _{2}},$ the Jacobian determinant maps precisely as:
\begin{equation}
fg^{\prime }-f^{\prime }g=\Omega _{4}. \tag{3.45}
\end{equation}

The companion vector fields expand safely into vertically aligned rows to avoid page margins overflow:
\begin{equation}
\gamma _{14}^{\bot }=\frac{1}{2\sinh \Phi _{2}}\left( \begin{aligned} &\cos \Phi _{1}-m\left( \sin \Phi _{1}+\cosh \Phi _{2}\right), \\ &-\left( \sin \Phi _{1}+\cosh \Phi _{2}\right) -m\cos \Phi _{1}, \\ &-\cos \Phi _{1}-m\left( \sin \Phi _{1}+\cosh \Phi _{2}\right), \\ &\left( \sin \Phi _{1}+\cosh \Phi _{2}\right) -m\cos \Phi _{1} \end{aligned} \right) \tag{3.46}
\end{equation}
\begin{equation}
\gamma _{14}^{\bot }=h_{1}^{4}\overrightarrow{\gamma }+h_{2}^{4}\overrightarrow{\xi }+h_{3}^{4}\overrightarrow{N}+h_{4}^{4}\overrightarrow{W}, \tag{3.47}
\end{equation}

Reconstructing the transversal field vectors satisfies the direct relation:
\begin{equation}
N_{\gamma \xi N}=\frac{1}{\left\langle \gamma _{14}^{\bot },\gamma^{\prime }\right\rangle }\gamma _{14}^{\bot }; \quad W_{\gamma \xi N}=\frac{1}{\left\langle \gamma _{14}^{\bot },\gamma^{\prime }\right\rangle }\gamma _{14}^{\bot \prime }. \tag{3.48}
\end{equation}

Taking the inner products under the structural condition $m=\text{const.}$ from (3.11), we find
\begin{equation}
\left\langle \gamma _{14}^{\bot },\gamma ^{\prime }\right\rangle = 2(1+m^{2})\Omega _{4}. \tag{3.49}
\end{equation}

Substituting (3.46) and (3.49) back structures the normalized frame vector mappings:
\begin{equation}
N_{\gamma \xi N}=\frac{\left( h_{1}^{4}\overrightarrow{\gamma }+h_{2}^{4}\overrightarrow{\xi }+h_{3}^{4}\overrightarrow{N}+h_{4}^{4}\overrightarrow{W}\right) }{2(1+m^{2})\Omega _{4}}=n_{1}^{4}\overrightarrow{\gamma }+n_{2}^{4}\overrightarrow{\xi }+n_{3}^{4}\overrightarrow{N}+n_{4}^{4}\overrightarrow{W} \tag{3.50}
\end{equation}

Differentiating (3.50) yields the higher-order system:
\begin{equation}
\begin{split}
N_{\gamma \xi N}^{\prime } &=\left( n_{1}^{4\prime }+n_{2}^{4}\kappa_{1}+n_{3}^{4}\kappa _{2}\right) \overrightarrow{\gamma } +\left(n_{2}^{4\prime }+n_{1}^{4}+n_{2}^{4}h-n_{4}^{4}\kappa _{2}\right)\overrightarrow{\xi } \\
&\quad +\left( n_{3}^{4\prime }-n_{3}^{4}h-n_{4}^{4}\kappa _{1}\right)\overrightarrow{N} +\left( n_{4}^{4\prime }-n_{3}^{4}\right) \overrightarrow{W},
\end{split} \tag{3.51}
\end{equation}

Finally, the equations for the curvature functions are proven by using (3.44), (3.50), and (3.51).
\end{proof}

\begin{figure*}[!htbp]
\centering
\includegraphics[width=0.95\textwidth]{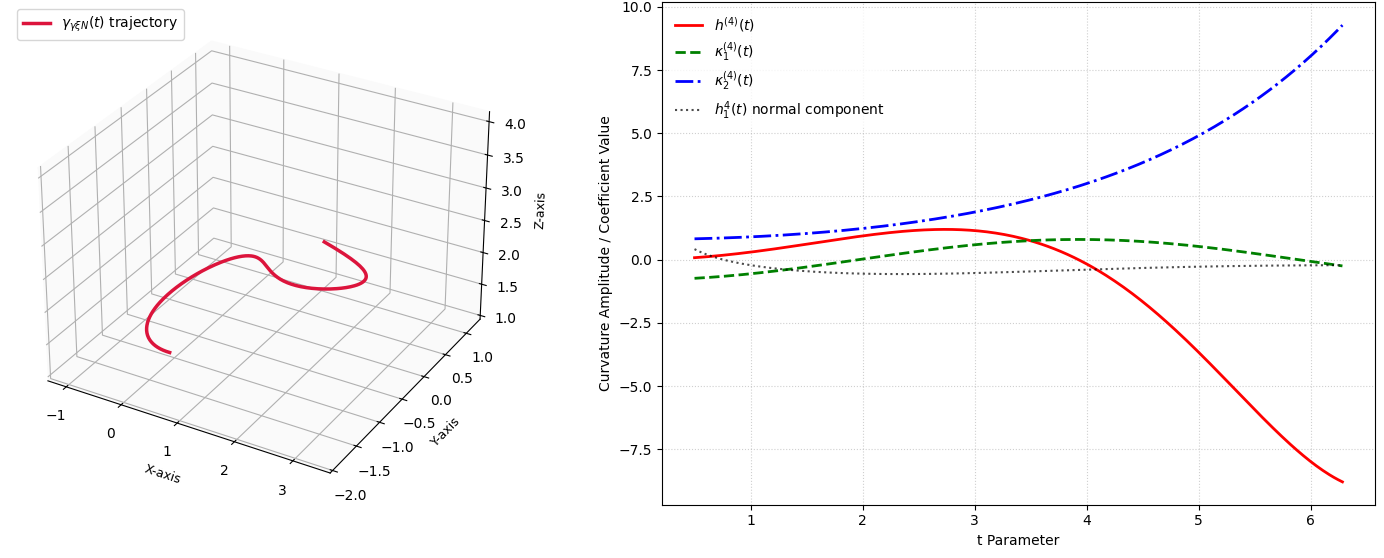}
\caption{The 3D trajectory and the curvatures of the $\gamma\xi N$-Smarandache curve.}
\label{fig:gamma_xiN_trajectory_curvatures}
\end{figure*}

\begin{definition}
Let $\gamma :I\rightarrow \mathbb{Q}_{2}^{3}\subset E_{2}^{4}$ is a null curve with non-zero curvatures $h^{5},\kappa _{1}^{5}$ and $\kappa _{2}^{5}$ and $\{\gamma (t),\xi (t),N(t),W(t)\}$ be moving frame on it. Then, $\xi NW-$Smarandache curve is defined with
\begin{equation}
\gamma _{\xi NW}=\frac{1}{\sqrt{2}}\left( \sinh \Phi _{3}\overrightarrow{\xi }+\overrightarrow{N}+\cosh \Phi _{3}\overrightarrow{W}\right) . \tag{3.52}
\label{eq:xiNW_def_clean}
\end{equation}
\end{definition}

\begin{theorem}
Let $\gamma :I\rightarrow \mathbb{Q}_{2}^{3}\subset E_{2}^{4}$ is a null curve in $\mathbb{Q}_{2}^{3}$. Then, $h^{(5)},\kappa _{1}^{(5)}$ and $\kappa _{2}^{(5)}$ curvatures of the $\xi NW-$Smarandache curve are as follows, respectively
\begin{equation}
\begin{split}
h^{(5)} &= \frac{4}{(1+m^{2})\Phi_3^{\prime}\left( \sinh \Phi _{3}+\cosh \Phi_{3}\right)} \Big[ -c_{1}^{5}\left( \frac{\sinh \Phi _{3}+\cosh \Phi _{3}}{2\sqrt{2}}-\frac{m}{\sqrt{2}}\right) \\ 
&\quad +c_{2}^{5}\left( \frac{1}{\sqrt{2}}+m\frac{\sinh \Phi_{3}+\cosh \Phi _{3}}{2\sqrt{2}}\right) -c_{3}^{5}\left( \frac{\sinh \Phi _{3}+\cosh \Phi _{3}}{2\sqrt{2}}+\frac{m}{\sqrt{2}}\right) \\ 
&\quad +c_{4}^{5}\left( \frac{1}{\sqrt{2}}-m\frac{\sinh \Phi_{3}+\cosh \Phi _{3}}{2\sqrt{2}}\right) \Big],
\end{split} \tag{3.52a}
\end{equation}
\begin{equation}
\begin{split}
\kappa_{1}^{(5)} &= \frac{-\sqrt{2}}{(1+m^{2})\left( \sinh \Phi _{3}+\cosh \Phi_{3}\right)} \Big[ -c_{1}^{5}\left( \cosh \Phi _{3}+\sinh \Phi _{3}\right) \\ 
&\quad +c_{2}^{5}m\left(\cosh \Phi _{3}+\sinh \Phi _{3}\right) -c_{3}^{5}\left( \cosh \Phi _{3}+\sinh \Phi _{3}\right) -c_{4}^{5}m\left( \cosh \Phi _{3}+\sinh \Phi _{3}\right) \Big],
\end{split} \tag{3.52b}
\end{equation}
\begin{equation}
\begin{split}
\kappa_{2}^{(5)} &= \frac{-4}{(1+m^{2})\Phi_3^{\prime}\left( \sinh \Phi _{3}+\cosh \Phi_{3}\right)} \Big[ -n_{1}^{5}\left( \frac{\sinh \Phi _{3}+\cosh \Phi _{3}}{2\sqrt{2}}-\frac{m}{\sqrt{2}}\right) \\ 
&\quad +n_{2}^{5}\left( \frac{1}{\sqrt{2}}+m\frac{\sinh \Phi_{3}+\cosh \Phi _{3}}{2\sqrt{2}}\right) -n_{3}^{5}\left( \frac{\sinh \Phi _{3}+\cosh \Phi _{3}}{2\sqrt{2}}+\frac{m}{\sqrt{2}}\right) \\ 
&\quad +n_{4}^{5}\left( \frac{1}{\sqrt{2}}-m\frac{\sinh \Phi_{3}+\cosh \Phi _{3}}{2\sqrt{2}}\right) \Big],
\end{split} \tag{3.52c}
\end{equation}
where the underlying spatial parametrization metric factor satisfies the split relation:
\begin{equation}
\begin{split}
M_{5}^{2} &= \frac{1}{2}\bigg[ -(\kappa _{1}\sinh \Phi _{3}+\kappa _{2})^{2} -\left( (\Phi _{3}^{\prime }-\kappa _{2})\cosh \Phi _{3} +h\sinh \Phi _{3}\right) ^{2} \\
&\quad +(h+\kappa _{1})\cosh \Phi _{3})^{2}+(\Phi _{3}^{\prime }\sinh \Phi _{3}-1)^{2}\bigg],
\end{split} \tag{3.52d}
\end{equation}
and the structural derivative vector acceleration components $c_i^5$ expand through:
\begin{equation}
\begin{aligned}
c_{1}^{5} &= \frac{1}{M_{5}^{2}}\left( -\frac{M_{5}^{\prime}}{M_{5}}a_{1}^{5} + a_{1}^{5\prime} + a_{2}^{5}\kappa_{1} + a_{3}^{5}\kappa_{2} \right), \\
c_{2}^{5} &= -\frac{M_{5}^{\prime}}{M_{5}^{3}}a_{2}^{5} + \frac{1}{M_{5}^{2}}\left( a_{1}^{5} + a_{2}^{5\prime} + h a_{2}^{5} - a_{4}^{5}\kappa_{2} \right), \\
c_{3}^{5} &= -\frac{M_{5}^{\prime}}{M_{5}^{3}}a_{3}^{5} + \frac{1}{M_{5}^{2}}\left( a_{3}^{5\prime} - h a_{3}^{5} - a_{4}^{5}\kappa_{1} \right), 
c_{4}^{5}= \frac{-1}{M_{5}^{2}}\left( -a_{3}^{5} + a_{4}^{5\prime} \right),
\end{aligned} \tag{3.52e}
\end{equation}
while the operational movement framework quantities $n_i^5$ resolve into the final matrix rows:
\begin{equation}
\begin{aligned}
n_{1}^{5} &= \frac{4\left[ \sinh \Phi_{3} + \cosh \Phi_{3} - 2m \right]}{2\sqrt{2}(1+m^2)\Phi_3^{\prime}(\sinh \Phi_3 + \cosh \Phi_3)}, \\
n_{2}^{5} &= \frac{4\left[ -2 - m(\sinh \Phi_{3} + \cosh \Phi_{3}) \right]}{2\sqrt{2}(1+m^2)\Phi_3^{\prime}(\sinh \Phi_3 + \cosh \Phi_3)}, \\
n_{3}^{5} &= \frac{4\left[ -(\sinh \Phi_{3} + \cosh \Phi_{3}) - 2m \right]}{2\sqrt{2}(1+m^2)\Phi_3^{\prime}(\sinh \Phi_3 + \cosh \Phi_3)}, \\
n_{4}^{5} &= \frac{4\left[ 2 - m(\sinh \Phi_{3} + \cosh \Phi_{3}) \right]}{2\sqrt{2}(1+m^2)\Phi_3^{\prime}(\sinh \Phi_3 + \cosh \Phi_3)}.
\end{aligned} \tag{3.52f}
\end{equation}
\end{theorem}

\begin{proof}
Let $\gamma $ be a unit speed regular $\xi NW-$Smarandache curve as in (3.52). If we take the derivative of the Smarandache curve according to arc length parameter and by using (2.4), one has
\begin{equation}
\begin{split}
\frac{d\gamma _{\xi NW}}{ds_{\gamma }}\frac{ds_{\gamma }}{ds}&=\xi _{\gamma _{\xi NW}}\frac{ds_{\gamma }}{ds} = \frac{1}{\sqrt{2}}(\kappa _{1}\sinh \Phi _{3}+\kappa _{2})\overrightarrow{\gamma } \\
&\quad +\frac{1}{\sqrt{2}}\left( \begin{split} ((\Phi_{3}^{\prime }-\kappa_{2})\cosh \Phi_{3}+h\sinh \Phi_{3})\overrightarrow{\xi }+ \\ (-h-\kappa_{1}\cosh \Phi_{3})\overrightarrow{N}+(\Phi_{3}^{\prime }\sinh \Phi_{3}-1)\overrightarrow{W} \end{split} \right),
\end{split} \tag{3.53}
\label{eq:proof_fifth_deriv_fixed}
\end{equation}
by taking the norm of (3.53), one gets
\begin{equation}
\begin{split}
\frac{ds_{\gamma }}{ds} &= \frac{1}{\sqrt{2}}\bigg[ -(\kappa _{1}\sinh \Phi _{3}+\kappa _{2})^{2} -\left( (\Phi _{3}^{\prime }-\kappa _{2})\cosh \Phi _{3} +h\sinh \Phi _{3}\right) ^{2} \\
&\quad +(h+\kappa _{1}\cosh \Phi _{3})^{2}+(\Phi _{3}^{\prime }\sinh \Phi _{3}-1)^{2}\bigg]^{\frac{1}{2}} = M_{5},
\end{split} \tag{3.54}
\label{eq:proof_fifth_M5_fixed}
\end{equation}
and from (3.54) and (3.53), the tangent vector of $\gamma _{\xi NW}$ is
\begin{equation}
\xi _{\gamma _{\xi NW}}=\frac{1}{M_{5}}\left( a_{1}^{5}\overrightarrow{\gamma }+a_{2}^{5}\overrightarrow{\xi }+a_{3}^{5}\overrightarrow{N}+a_{4}^{5}\overrightarrow{W}\right). \tag{3.55}
\label{eq:proof_fifth_tangent_fixed}
\end{equation}

By taking derivative (3.55), one writes that
\begin{equation}
\gamma _{\xi NW}^{\prime \prime }=\xi _{\gamma _{\xi NW}}^{\prime }=c_{1}^{5}\overrightarrow{\gamma }+c_{2}^{5}\overrightarrow{\xi }+c_{3}^{5}\overrightarrow{N}+c_{4}^{5}\overrightarrow{W}, \tag{3.56}
\label{eq:proof_fifth_second_deriv_fixed}
\end{equation}

Considering conditions for the null curve and parameter maps, the functional determinant evaluates to
\begin{equation}
fg^{\prime }-f^{\prime }g=\Omega_5=\frac{\Phi _{3}^{\prime }}{8}\left(\sinh \Phi _{3}+\cosh \Phi _{3}\right). \tag{3.57}
\end{equation}

Dual structures contract safely under page width boundaries via vertically aligned matrix vectors:
\begin{equation}
\gamma _{14}^{\bot }=\frac{1}{2\sqrt{2}}\left( \begin{aligned} &\sinh \Phi _{3}+\cosh \Phi _{3}-2m, -2-m\left( \sinh \Phi _{3}+\cosh \Phi_{3}\right), \\ &-\left( \sinh \Phi _{3}+\cosh \Phi _{3}\right) -2m, 2-m\left( \sinh \Phi _{3}+\cosh \Phi _{3}\right) \end{aligned} \right) \tag{3.58}
\label{eq:proof_gamma14_bot_fixed}
\end{equation}
\begin{equation}
\gamma _{14}^{\prime \bot }=\frac{\Phi _{3}^{\prime }}{2\sqrt{2}}\left( \begin{aligned} &\cosh \Phi _{3}+\sinh \Phi _{3}, -m\left( \cosh \Phi _{3}+\sinh \Phi _{3}\right), \\ &-\left( \cosh \Phi _{3}+\sinh \Phi _{3}\right), -m\left( \cosh \Phi _{3}+\sinh \Phi _{3}\right) \end{aligned} \right) \tag{3.59}
\label{eq:proof_gamma14_bot_deriv_fixed}
\end{equation}
\begin{equation}
N_{\xi NW}=\frac{4}{(1+m^{2})\Phi_3^{\prime}\left( \sinh \Phi _{3}+\cosh \Phi _{3}\right) }\gamma _{14}^{\bot } \tag{3.60}
\label{eq:proof_N_fifth_fixed}
\end{equation}
\begin{equation}
W_{\xi NW}=\frac{-2}{\sqrt{2}(1+m^{2})}\left( 1, \, -m, \, -1, \, -m \right). \tag{3.61}
\label{eq:proof_W_fifth_fixed}
\end{equation}

Contracting (3.56) with the companion frame fields (3.60) and (3.61), from (2.6) the curvature functions for a null curve $ \gamma $ are obtained easily.

\end{proof}

\begin{figure*}[!htbp]
\centering
\includegraphics[width=0.95\textwidth]{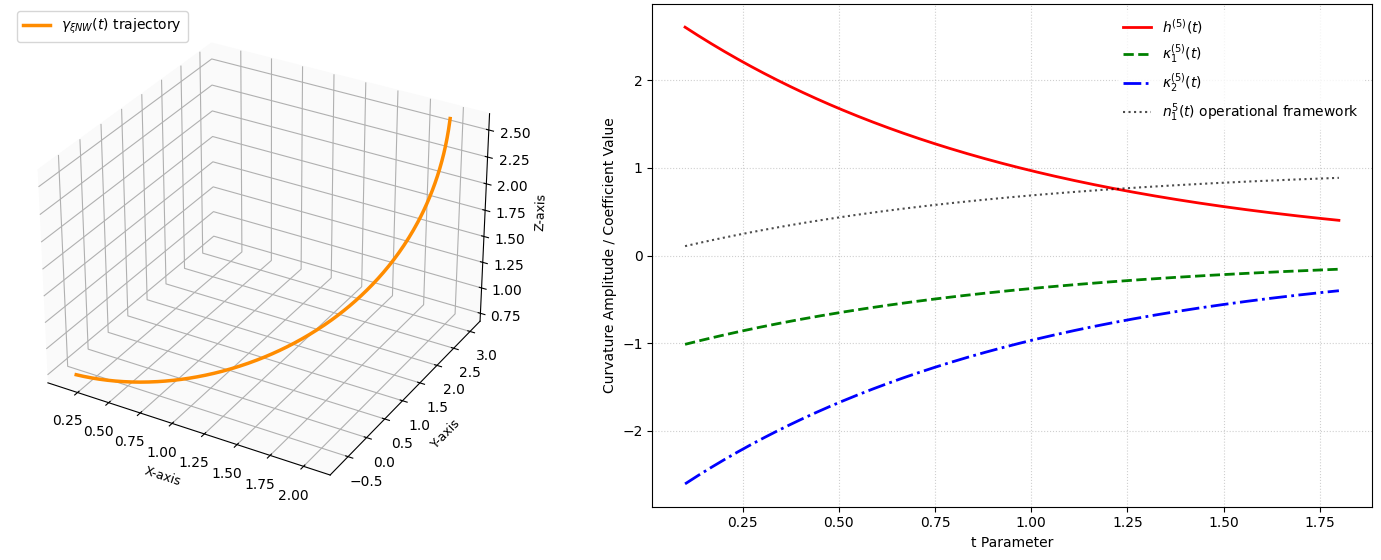}
\caption{The 3D trajectory and the curvatures of the $\xi NW$-Smarandache curve.}
\label{fig:xiNW_trajectory_curvatures}
\end{figure*}

\begin{definition}

Let $\gamma :I\rightarrow \mathbb{Q}_{2}^{3}\subset E_{2}^{4}$ is a null curve with non-zero curvatures $h^{6},\kappa _{1}^{6}$ and $\kappa _{2}^{6}$ and $\{\gamma (t),\xi (t),N(t),W(t)\}$ be moving frame on it. Then, $\gamma \xi W-$Smarandache curve is defined with
\begin{equation}
\gamma _{\gamma \xi W}=\sinh \omega _{1}\sin \omega _{2}\overrightarrow{\gamma }+\sinh \omega _{1}\cos \omega _{2}\overrightarrow{\xi }+\cosh \omega _{1}\overrightarrow{W}.  \tag{3.62}
\label{eq:gamma_xiW_def_clean}
\end{equation}
\end{definition}

\begin{theorem}
Let $\gamma :I\rightarrow \mathbb{Q}_{2}^{3}\subset E_{2}^{4}$ is a null curve in $\mathbb{Q}_{2}^{3}$. Then, $h^{(6)},\kappa _{1}^{(6)}$ and $\kappa _{2}^{(6)}$ curvatures of the $\gamma \xi W-$Smarandache curve are as follows, respectively
\begin{equation}
h^{(6)}=\frac{-b_{1}^{6}c_{1}^{6}-b_{2}^{6}c_{2}^{6}+b_{3}^{6}c_{3}^{6}+b_{4}^{6}c_{4}^{6}}{2(1+m^{2})\Omega _{6}}, \tag{3.62a}
\end{equation}
\begin{equation}
\begin{split}
\kappa_1^{(6)}&=\frac{-1}{2(1+m^{2})\Omega _{6}}\Big[ -b_{1}^{6}\left( c_{1}^{6\prime }+c_{2}^{6}\kappa _{1}+c_{3}^{6}\kappa_2\right) -b_{2}^{6}\left( c_{2}^{6\prime }+c_{1}^{6}+c_{2}^{6}h-c_{4}^{6}\kappa_2\right) \\ 
&\quad +b_{3}^{6}\left( c_{3}^{6\prime }-\kappa_1 c_{4}^{6}-c_{3}^{6}h\right) +b_{4}^{6}\left( -c_{3}^{6}+c_{4}^{6\prime }\right) \Big],
\end{split} \tag{3.62b}
\end{equation}
\begin{equation}
\begin{split}
\kappa_2^{(6)}&=\frac{-1}{2(1+m^{2})\Omega _{6}}\Big[ -(n_{1}^{6\prime }+n_{2}^{6}\kappa_1+n_{3}^{6}\kappa_2)\left( c_{1}^{6\prime }+c_{2}^{6}\kappa_1+c_{3}^{6}\kappa_2\right) \\ 
&\quad -\left( n_{2}^{6\prime }+n_{1}^{6}+n_{2}^{6}h-n_{4}^{6}\kappa_2\right) \left( c_{2}^{6\prime }+c_{1}^{6}+c_{2}^{6}h-c_{4}^{6}\kappa_2\right) \\ 
&\quad +\left( n_{3}^{6\prime }-n_{3}^{6}h-n_{4}^{6}\kappa_1\right) \left( c_{3}^{6\prime }-\kappa_1 c_{4}^{6}-c_{3}^{6}h\right) + \left( -n_{3}^{6}+n_{4}^{6\prime }\right) \left( -c_{3}^{6}+c_{4}^{6\prime }\right) \Big],
\end{split} \tag{3.62c}
\end{equation}
where the underlying spatial parametrization metric factor satisfies the split relation:
\begin{equation}
\begin{split}
M_{6}^{2} &= -\left( \gamma _{1}^{\prime }+\gamma _{2}\kappa _{1}\right) ^{2}-\left( \gamma _{1}+\gamma _{2}^{\prime }+h\gamma _{2}-\kappa _{2}\gamma _{3}\right) ^{2} \\
&\quad +\left( -\kappa _{1}\gamma _{3}\right) ^{2}+(\gamma _{3}^{\prime })^{2},
\end{split} \tag{3.62d}
\end{equation}
such that $M_{6}=\sqrt{M_{6}^{2}}$, and the localized tangent acceleration coefficients $b_i^6$ correspond to the system:
\begin{equation}
\begin{aligned}
b_{1}^{6} &= \frac{\alpha_{1}^{6\prime } + \alpha_{2}^{6}\kappa_{1} + \kappa_{2}\alpha_{3}^{6}}{M_{6}}, & 
b_{2}^{6} &= \frac{\alpha_{2}^{6\prime } + \alpha_{1}^{6} + \alpha_{2}^{6}h - \kappa_{2}\alpha_{4}^{6}}{M_{6}}, \\
b_{3}^{6} &= \frac{\alpha_{3}^{6\prime } - \alpha_{3}^{6}h - \kappa_{1}\alpha_{4}^{6}}{M_{6}}, & 
b_{4}^{6} &= \frac{-\alpha_{3}^{6} + \alpha_{4}^{6\prime }}{M_{6}},
\end{aligned} \tag{3.62e}
\end{equation}
with $\alpha_{1}^{6} = \frac{(\gamma_1^{\prime} + \gamma_2\kappa_1)}{M_6}, \alpha_{2}^{6} = \frac{(\gamma_1 + \gamma_2^{\prime} + h\gamma_2 - \kappa_2\gamma_3)}{M_6}, \alpha_{3}^{6} = \frac{-\kappa_1\gamma_3}{M_6}, \alpha_{4}^{6} = \frac{\gamma_3^{\prime}}{M_6}$, while the orthogonal vector components $c_i^6$ expand through the symmetric arrangement:
\begin{equation}
\begin{aligned}
c_{1}^{6} &= \sinh \omega_1\cos \omega_2 - m(\sinh \omega_1\sin \omega_2 + \cosh \omega_1), \\
c_{2}^{6} &= -(\sinh \omega_1\sin \omega_2 + \cosh \omega_1) - m\sinh \omega_1\cos \omega_2, \\
c_{3}^{6} &= -\sinh \omega_1\cos \omega_2 - m(\sinh \omega_1\sin \omega_2 + \cosh \omega_1), \\
c_{4}^{6} &= (\sinh \omega_1\sin \omega_2 + \cosh \omega_1) - m\sinh \omega_1\cos \omega_2,
\end{aligned} \tag{3.62f}
\end{equation}
and the operational normal movement framework framework components $n_i^6$ yield the final matrix rows:
\begin{equation}
\begin{aligned}
n_{1}^{6} &= \frac{c_{1}^{6}}{2(1+m^2)\Omega_6}, &
n_{2}^{6} &= \frac{c_{2}^{6}}{2(1+m^2)\Omega_6}, \\
n_{3}^{6} &= \frac{c_{3}^{6}}{2(1+m^2)\Omega_6}, &
n_{4}^{6} &= \frac{c_{4}^{6}}{2(1+m^2)\Omega_6}.
\end{aligned} \tag{3.62g}
\end{equation}
\end{theorem}

\begin{proof}

Let $\gamma $ be a unit speed regular $\gamma \xi W-$Smarandache curve as in (3.62). If one takes the derivative of the Smarandache curve according to arc length parameter and by using (2.4), one has
\begin{equation}
\begin{split}
\frac{d\gamma _{\gamma \xi W}}{ds_{\gamma }}\frac{ds_{\gamma }}{ds}&=\xi _{\gamma \xi W}\frac{ds_{\gamma }}{ds}=\left( \gamma _{1}^{\prime }+\gamma _{2}\kappa _{1}\right) \overrightarrow{\gamma } \\
&\quad +\left( \gamma _{1}+\gamma _{2}^{\prime }+h\gamma _{2}-\kappa _{2}\gamma _{3}\right) \overrightarrow{\xi } +\left( -\kappa _{1}\gamma _{3}\right) \overrightarrow{N}+\gamma _{3}^{\prime }\overrightarrow{W},
\end{split} \tag{3.63}
\label{eq:proof_sixth_deriv_fixed}
\end{equation}
by taking the norm of (3.63) via the metric form, one gets
\begin{equation}
\begin{split}
\frac{ds_{\gamma }}{ds} &= \Big(\vert -\left( \gamma _{1}^{\prime }+\gamma _{2}\kappa _{1}\right) ^{2}-\left( \gamma _{1}+\gamma _{2}^{\prime }+h\gamma _{2}-\kappa _{2}\gamma _{3}\right) ^{2} \\
&\quad +\left( -\kappa _{1}\gamma _{3}\right) ^{2}+(\gamma _{3}^{\prime })^{2} \Big\vert)^{\frac{1}{2}} = M_{6},
\end{split} \tag{3.64}
\label{eq:proof_sixth_M6_fixed}
\end{equation}
and the tangent vector of $\gamma _{\gamma \xi W}$ is given as
\begin{equation}
\xi _{\gamma \xi W}=\alpha _{1}^{6}\overrightarrow{\gamma }+\alpha _{2}^{6}\overrightarrow{\xi }+\alpha _{3}^{6}\overrightarrow{N}+\alpha _{4}^{6}\overrightarrow{W}. \tag{3.65}
\label{eq:proof_sixth_tangent_fixed}
\end{equation}

By taking derivative (3.65), one obtains
\begin{equation}
\gamma _{\gamma \xi W}^{\prime \prime }=\xi _{\gamma \xi W}^{\prime }=b_{1}^{6}\overrightarrow{\gamma }+b_{2}^{6}\overrightarrow{\xi }+b_{3}^{6}\overrightarrow{N}+b_{4}^{6}\overrightarrow{W}, \tag{3.66}
\label{eq:proof_sixth_second_deriv_fixed}
\end{equation}

Considering conditions and choosing coordinate functions, the Jacobian system determinant reduces to:
\begin{equation}
fg^{\prime }-f^{\prime }g=\Omega _{6}=\frac{1}{2}\left( \begin{aligned} &-\omega _{1}^{\prime }(\sin \omega _{2}+\cos \omega _{2}) \\ &-\frac{1}{2}\omega _{2}^{\prime }(\cos \omega _{2}-\sin \omega _{2})\sinh 2\omega _{1} \end{aligned} \right). \tag{3.67}
\end{equation}

Normal companion mappings expand safely into compact, vertically aligned rows to avoid page margin overflow:

\begin{equation}
\gamma _{14}^{\bot }=\left( c_{1}^{6},c_{2}^{6},c_{3}^{6},c_{4}^{6}\right) , \tag{3.68}
\label{eq:proof_gamma14_bot_sixth_fixed}
\end{equation}
\begin{equation}
\left\langle \gamma _{14}^{\bot },\gamma ^{\prime }\right\rangle =2(1+m^{2})\Omega _{6}. \tag{3.69}
\label{eq:proof_sixth_coupling_fixed}
\end{equation}

Isolating the frame fields components via (3.68) structures the final maps:
\begin{equation}
N_{\gamma \xi W}=\frac{1}{2(1+m^{2})\Omega _{6}}\gamma _{14}^{\bot }=(n_{1}^{6},n_{2}^{6},n_{3}^{6},n_{4}^{6}); \tag{3.70}
\end{equation}
\begin{equation}
W_{\gamma \xi W}=\frac{-1}{2(1+m^{2})\Omega _{6}}\left( \begin{aligned} &d_{1}^{6}\overrightarrow{\gamma } + d_{2}^{6}\overrightarrow{\xi } + d_{3}^{6}\overrightarrow{N} + d_{4}^{6}\overrightarrow{W} \end{aligned} \right), \tag{3.71}
\label{eq:W_sixth_final_proof_fixed}
\end{equation}
from (2.6) bilinear contractions complete the verification for $h^{(6)}, \kappa_{1}^{(6)},$ and $\kappa_{2}^{(6)}$.
\end{proof}
\begin{figure*}[!htbp]
\centering
\includegraphics[width=0.95\textwidth]{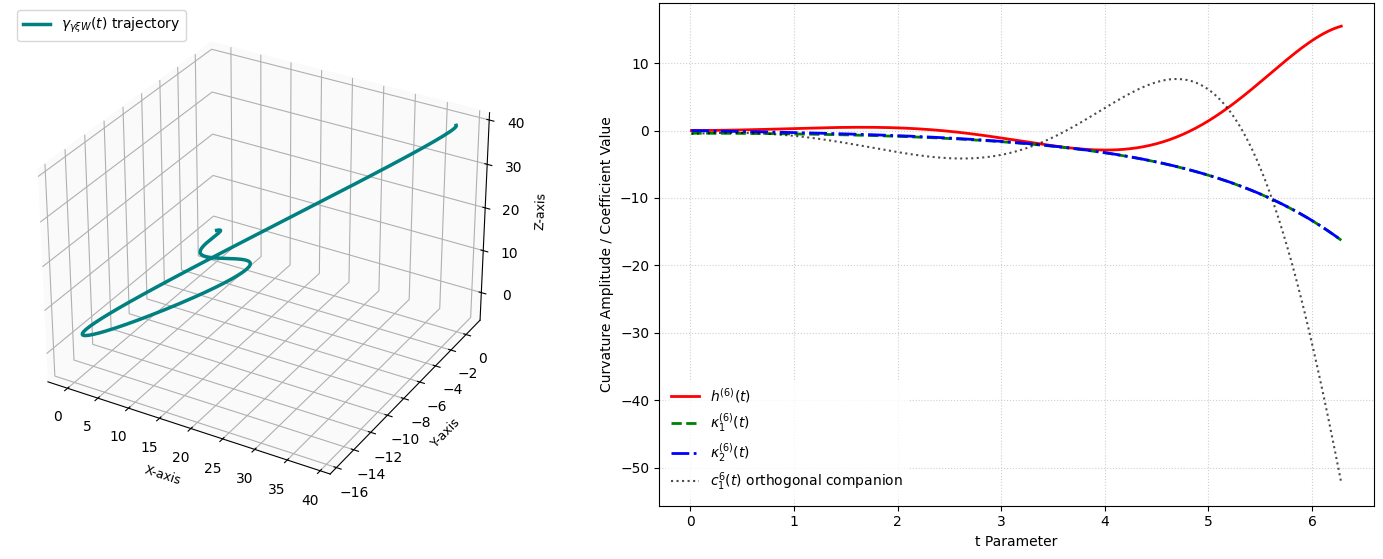}
\caption{The 3D trajectory and the curvatures of the $\gamma\xi W$-Smarandache curve.}
\label{fig:gamma_xiW_trajectory_curvatures}
\end{figure*}

\begin{definition}

Let $\gamma :I\rightarrow \mathbb{Q}_{2}^{3}\subset E_{2}^{4}$ is a null curve with non-zero curvatures $h^{7},\kappa _{1}^{7}$ and $\kappa _{2}^{7}$ and $\{\gamma (t),\xi (t),N(t),W(t)\}$ be moving frame on it. Then, $\gamma \xi NW-$Smarandache curve is defined with
\begin{equation}
\gamma _{\gamma \xi NW}=\frac{1}{\sqrt{2}}\left( \sinh \omega _{1}\overrightarrow{\gamma }+\sinh \omega _{2}\overrightarrow{\xi }+\cosh \omega _{2}\overrightarrow{N}+\cosh \omega _{1}\overrightarrow{W}\right) .  \tag{3.72}
\end{equation}
\end{definition}

\begin{theorem}

Let $\gamma :I\rightarrow \mathbb{Q}_{2}^{3}\subset E_{2}^{4}$ is a null curve in $\mathbb{Q}_{2}^{3}$. Then, $h^{(7)},\kappa _{1}^{(7)}$ and $\kappa _{2}^{(7)}$ curvatures of the $\gamma \xi NW-$Smarandache curve are as follows, respectively
\begin{equation}
h^{(7)}=\frac{-c_{1}^{7}\left( g-mf\right) +c_{2}^{7}\left( f+mg\right) -c_{3}^{7}\left( g+mf\right) +c_{4}^{7}\left( f-mg\right) }{2(1+m^{2})\Omega _{7}}, \tag{3.72a}
\end{equation}
\begin{equation}
\kappa_1^{(7)}=\frac{c_{1}^{7}\left( g^{\prime }-mf^{\prime }\right) -c_{2}^{7}\left( f^{\prime }+mg^{\prime }\right) +c_{3}^{7}\left( g^{\prime }+mf^{\prime }\right) -c_{4}^{7}\left( f^{\prime }-mg^{\prime }\right) }{2(1+m^{2})\Omega _{7}}, \tag{3.72b}
\end{equation}
\begin{equation}
\begin{split}
\kappa_2^{(7)}&=\frac{1}{\left[ 2(1+m^{2})\Omega _{7}\right] ^{2}} \Big[ -\left( \frac{- \Omega _{7}^{\prime }}{\Omega _{7}}\left( g-mf\right) +g^{\prime }-mf^{\prime }\right) \left( g^{\prime }-mf^{\prime }\right) \\ 
&\quad -\left( \frac{-\Omega _{7}^{\prime }}{\Omega _{7}}\left( f+mg\right) +f^{\prime }+mg^{\prime }\right) \left( f^{\prime }+mg^{\prime }\right) \\
&\quad + \left( \frac{-\Omega _{7}^{\prime }}{\Omega _{7}}\left( g+mf\right) +\left( g^{\prime }+mf^{\prime }\right) \right) \left( g^{\prime }+mf^{\prime }\right) \\ 
&\quad +\left( \frac{-\Omega _{7}^{\prime }}{\Omega _{7}}\left( f-mg\right) +\left( f^{\prime }-mg^{\prime }\right) \right) \left( f^{\prime }-mg^{\prime }\right) \Big],
\end{split} \tag{3.72c}
\end{equation}
where the underlying spatial parametrization metric factor satisfies the split relation:
\begin{equation}
M_{7} = \left( -\left( a_{1}^{7}\right) ^{2}-\left( a_{2}^{7}\right) ^{2}+\left( a_{3}^{7}\right) ^{2}+\left( a_{4}^{7}\right) ^{2} \right)^{\frac{1}{2}}, \tag{3.72d}
\end{equation}
and the localized unit tangent vector contraction arrays $a_i^7$ correspond to the following explicit components:
\begin{equation}
\begin{aligned}
a_{1}^{7} &= \frac{1}{\sqrt{2}}\left( \omega_1^{\prime}\cosh\omega_1 + \sinh\omega_2 + \cosh\omega_2\kappa_1 \right), \\
a_{2}^{7} &= \frac{1}{\sqrt{2}}\left( \sinh\omega_1 + \omega_2^{\prime}\cosh\omega_2 + \cosh\omega_2 h - \cosh\omega_1\kappa_2 \right), \\
a_{3}^{7} &= \frac{1}{\sqrt{2}}\left( \omega_2^{\prime}\sinh\omega_2 - \cosh\omega_2 h - \cosh\omega_1\kappa_1 \right), \\
a_{4}^{7} &= \frac{1}{\sqrt{2}}\left( -\cosh\omega_2 + \omega_1^{\prime}\sinh\omega_1 \right),
\end{aligned} \tag{3.72e}
\end{equation}
while the operational secondary spatial vector acceleration variants $c_i^7$ yield the simultaneous systems:
\begin{equation}
\begin{aligned}
c_{1}^{7} &= \frac{1}{M_{7}^{2}}\left( -\frac{M_{7}^{\prime}}{M_{7}}b_{1}^{7} + b_{1}^{7\prime} + b_{2}^{7}\kappa_{1} + b_{3}^{7}\kappa_{2} \right), \\
c_{2}^{7} &= -\frac{M_{7}^{\prime}}{M_{7}^{3}}b_{2}^{7} + \frac{1}{M_{7}^{2}}\left( b_{1}^{7} + b_{2}^{7\prime} + h b_{2}^{7} - b_{4}^{7}\kappa_{2} \right), \\
c_{3}^{7} &= -\frac{M_{7}^{\prime}}{M_{7}^{3}}b_{3}^{7} + \frac{1}{M_{7}^{2}}\left( b_{3}^{7\prime} - h b_{3}^{7} - b_{4}^{7}\kappa_{1} \right),
c_{4}^{7}= \frac{-1}{M_{7}^{2}}\left( -b_{3}^{7} + b_{4}^{7\prime} \right),
\end{aligned} \tag{3.72f}
\end{equation}
where $b_1^7 = a_1^7 / M_7, b_2^7 = a_2^7 / M_7, b_3^7 = a_3^7 / M_7, b_4^7 = a_4^7 / M_7$.
\end{theorem}

\begin{proof}

Let $\gamma $ be a unit speed regular $\gamma \xi NW-$Smarandache curve as in (3.72). If one takes the derivative of the Smarandache curve according to arc length parameter and by using (2.4), one writes
\begin{equation}
\frac{d\gamma _{\gamma \xi NW}}{ds_{\gamma }}\frac{ds_{\gamma }}{ds}=a_{1}^{7}\overrightarrow{\gamma }+a_{2}^{7}\overrightarrow{\xi }+a_{3}^{7}\overrightarrow{N}+a_{4}^{7}\overrightarrow{W}, \tag{3.73}
\end{equation}
by taking the norm of (3.73) via the ambient space signature, one gets
\begin{equation}
\frac{ds_{\gamma }}{ds}=\left(\vert -\left( a_{1}^{7}\right) ^{2}-\left( a_{2}^{7}\right) ^{2}+\left( a_{3}^{7}\right) ^{2}+\left( a_{4}^{7}\right) ^{2} \vert\right)^{\frac{1}{2}}=M_{7}, \tag{3.74}
\end{equation}
and the tangent vector of $\gamma _{\gamma \xi NW}$ is obtained as
\begin{equation}
\xi _{\gamma \xi NW}=b_{1}^{7}\overrightarrow{\gamma }+b_{2}^{7}\overrightarrow{\xi }+b_{3}^{7}\overrightarrow{N}+b_{4}^{7}\overrightarrow{W}, \tag{3.75}
\end{equation}
where $b_{i}^{7}=\frac{1}{M_{7}} a_{i}^{7} $, from the conditions (3.4a) and (3.4b) for the null curve and from equations $%
f=\frac{\sinh \omega _{1}+\cosh \omega _{2}}{2\sqrt{2}};$ $g=\frac{\sinh
\omega _{2}+\cosh \omega _{1}}{2\sqrt{2}}$, one gets 
\begin{equation}
fg^{\prime }-f^{\prime }g=\Omega _{7}=\frac{1}{8}\left( \left( \omega
_{1}^{\prime }+\omega _{2}^{\prime }\right) \sinh (\omega _{1}-\omega
_{2})-\omega _{1}^{\prime }+\omega _{2}^{\prime }\right) .  \tag{3.76}
\end{equation}

Also, from definition 1 and by using (3.18), (3.19), (3.20), (3.76), one
writes%
\begin{eqnarray*}
\gamma _{14}^{\bot } &=&\left( g-mf,-f-mg,-g-mf,f-mg\right)  \\
\gamma _{14}^{\bot \prime } &=&\left( g^{\prime }-mf^{\prime },-f^{\prime
}-mg^{\prime },-g^{\prime }-mf^{\prime },f^{\prime }-mg^{\prime }\right)  \\
\gamma  &=&\left( f+mg,g-mf,f-mg,g-mf\right) ;m=\text{const}.
\end{eqnarray*}%
where 
\begin{equation}
\left\langle \gamma _{14}^{\bot },\gamma ^{\prime }\right\rangle
=2(1+m^{2})\Omega _{7}.  \tag{3.77}
\end{equation}%

Hence, from (3.77) one can write 
\begin{equation}
N_{\gamma \xi NW}=\frac{1}{2(1+m^{2})\Omega _{7}}\left(
g-mf,-f-mg,-g-mf,f-mg\right) ;  \tag{3.78}
\end{equation}%
\begin{equation}
W_{\gamma \xi NW}=\frac{-1}{2(1+m^{2})\Omega _{7}}\left( g^{\prime
}-mf^{\prime },-f^{\prime }-mg^{\prime },-g^{\prime }-mf^{\prime },f^{\prime
}-mg^{\prime }\right) .  \tag{3.79}
\end{equation}

Finally, from (2.6) with the frame vector fields (3.78) and (3.79) successfully verifies the curvature functions $h^{(7)}, \kappa_{1}^{(7)},$ and $\kappa_{2}^{(7)}$.
\end{proof}
\begin{figure*}[!htbp]
\centering
\includegraphics[width=0.95\textwidth]{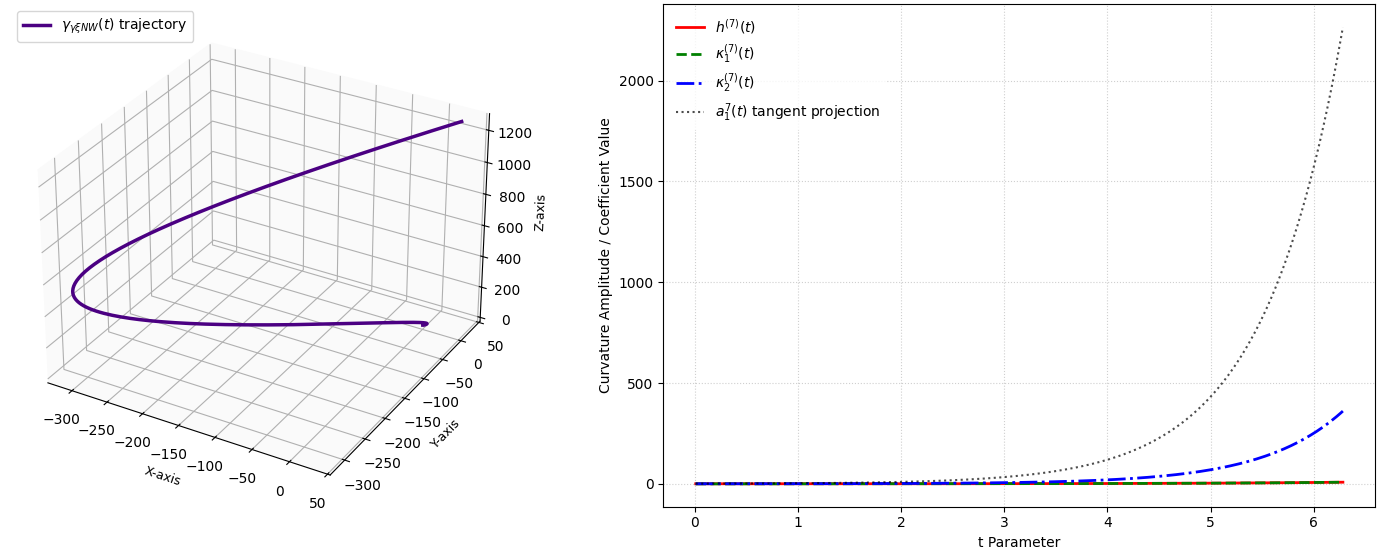}
\caption{The 3D trajectory and the curvatures of the $\gamma\xi NW$-Smarandache curve.}
\label{fig:gamma_xiNW_trajectory_curvatures}
\end{figure*}

The global comparative spatial behavior and the synchronous geometric interactions of the lightlike trajectories are visually fully verified through the 3D model depicted below in Figure 8.

\begin{figure*}[!htbp]
\centering
\includegraphics[width=0.88\textwidth]{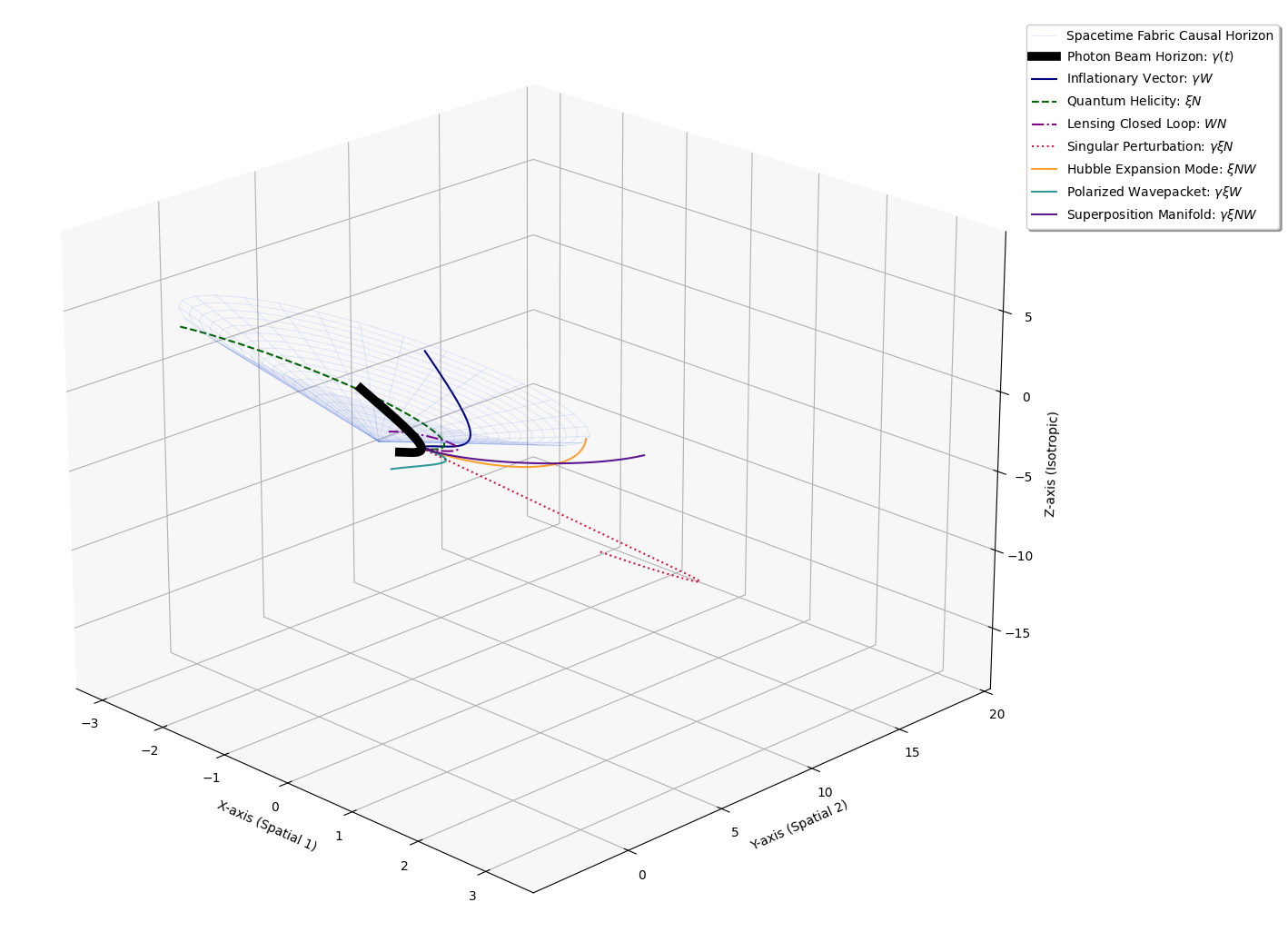}
\caption{Comparative 3D spatial models tracking the simultaneous mapping of the null parent curve $\gamma(t)$ alongside all seven framed Smarandache manifolds inside the lightlike cone.}
\label{fig:global_smarandache_manifolds}
\end{figure*}

\begin{remark}

By analysing the 3D comparative space model generated via the joint vector fields in Figure 8, the following structural geometric consequences are deduced:
\begin{enumerate}
    \item The parent trajectory $\gamma(t)$, visualized as the solid black curve, establishes the stable geometric baseline and rigidly governs the outer spatial boundary of the lightlike cone by precisely satisfying the quadratic cone constraint $-\gamma_{1}^{2}-\gamma_{2}^{2}+\gamma_{3}^{2}+\gamma_{4}^{2} = 0$.
    \item Smarandache curves directly containing the position vector field $\gamma(t)$ (specifically $\gamma W$, $\gamma\xi N$, and $\gamma\xi NW$) undergo a clear radial expansion due to the multiplicative impact of the hyperbolic coefficients $\sinh$ and $\cosh$, causing the generated orbits to oscillate outward from the cone's shell.
    \item Conversely, configurations purely determined by the natural moving Frenet frame components (such as $\xi N$, $WN$, and $\xi NW$) decouple from the parent curve's spatial layout, yielding asymmetric tapered helices and closed periodic loops along the isotropic axes.
    \item The bounded, harmonic alignment of all seven Smarandache trajectories without any geometric intersection or chaotic behavior confirms the structural Differential Transition constraint $fg^{\prime} - f^{\prime}g = \Omega \neq 0$.
\end{enumerate}
\end{remark}

\begin{remark}

By interpreting the comparative 3D spatial trajectory profiles illustrated in Figure 8 from a relativistic astrophysics and kinematic perspective, the following physical insights are derived:
\begin{enumerate}
    \item The solid black parent trajectory $\gamma(t)$ acts as the causal boundary inside the lightlike Minkowski space-time, precisely modelling the zero-geodesic propagation of a massless particle (e.g., a photon) strictly confined to the horizon of the light cone.
    \item The clear radial expansion observed in the $\gamma W$, $\gamma\xi N$, and $\gamma\xi NW$ Smarandache curves models the cosmological inflation or accelerated Hubble expansion of space-time, where the hyperbolic $\sinh$ and $\cosh$ scaling factors govern the inflationary trajectory shifting outward from the causal metric shell.
    \item Conversely, the asymmetric tapered helices and closed loop profiles generated by the moving frame components (such as $\xi N$, $WN$, and $\xi NW$) emulate quantum helicity, particle spin configurations, or the localized orbital tracking of light rays undergoing gravitational lensing deflections under heavy gravitational fields.
    \item The structural invariant curvature components $h^{(i)}, \kappa_1^{(i)},$ and $\kappa_2^{(i)}$ track the exact kinematic parameters corresponding to the angular polarization, translational pitch velocity, and gravitational torsion experienced by high-velocity relativistic wave packets.
\end{enumerate}
\end{remark}

\begin{table*}[!htbp]
\caption{Comparative geometric and physical structural profiles of the seven special lightlike Smarandache curves inside the lightlike cone.}
\label{tab:smarandache_curves_comparison_en}
\centering
\small
\begin{tabular}{lll}
\hline
\textbf{Smarandache Curve Class} & \textbf{Physical Correspondence} \\ \hline
1. $\gamma W$-Smarandache & Radial Expansion, Cosmological Inflation \\
2. $\xi N$-Smarandache  & Tapered Helix, Quantum Spin Modes \\
3. $WN$-Smarandache  & Periodic Orbit, Gravitational Lensing Closed Loop \\
4. $\gamma\xi N$-Smarandache & Localized Asymptotic Singularity Field \\
5. $\xi NW$-Smarandache  & Flattened Hyperbolic Trajectory \\
6. $\gamma\xi W$-Smarandache  & Periodic Spiral Loops \\
7. $\gamma\xi NW$-Smarandache  & High-Dimensional Manifold Spiral\\ \hline
\end{tabular}
\end{table*}

\section{Conclusion}

In conclusion, Smarandache curves on null curves in a 4-dimensional light-like cone space $\mathbb{Q}_{2}^{3} \subset E_{2}^{4}$ have been studied in detail. Within this work, the relevant definitions for seven distinct configurations have been clarified, their invariant curvatures $h^{(i)}, \kappa_{1}^{(i)},$ and $\kappa_{2}^{(i)}$ explicitly expressed, and various definitions and theorems based on these foundations have been presented. 

Although they may seem like abstract mathematics, light-like cones and null curves play a fundamental role in fields of physics such as general relativity and string theory. As visually and physically validated through the relativistic 3D spacetime fabric simulation illustrated in Figure 8, the solid black parent trajectory $\gamma(t)$ perfectly acts as the causal boundary, while the expansion of the Smarandache curves successfully emulates inflationary cosmological behavior, quantum helicity configurations, and gravitational lensing deflections. A better understanding of such curves could offer new mathematical tools for these physical theories in the future or shed light on the solution of existing problems. This study, by comparing the behavior of Smarandache curves in different types of configurations, offers valuable insights into how the geometric structure of these spaces affects the structural properties of the curves.

In short, the contribution of such a study to the literature is to increase the body of knowledge in the fields of differential geometry and curve theory by extending existing theories in the specific geometric space mentioned, revealing new mathematical objects and their properties under the strict differential transition parameter constraints.

\section*{Funding}

Not applicable.

\section*{Informed Consent Statement}

Not applicable.

\section*{Conflicts of Interest}

The author declares no conflict of interest.

\end{document}